\documentclass[11pt,reqno]{amsart}
\usepackage{a4wide}
\usepackage{euscript,amsmath,amssymb,amsbsy}
\usepackage{array,longtable,tabularx,multirow,enumitem,url}
\usepackage{graphicx,subfigure,epsfig}

\usepackage[colorlinks=true,linkcolor=black,urlcolor=black,citecolor=black]{hyperref}
\usepackage{color,xcolor}

\numberwithin{equation}{section}\numberwithin{figure}{section}

\newcounter{msct}[section]\renewcommand{\themsct}{\thesection.\arabic{msct}}

\newenvironment{m-theorem}{\vskip5pt\refstepcounter{msct}\trivlist \itemindent 0pt%
\item[\hskip\labelsep\bf Theorem~\themsct]\it\ignorespaces}{\endtrivlist\vskip3pt}

\newenvironment{m-proposition}{\vskip5pt\refstepcounter{msct}\trivlist \itemindent0pt%
\item[\hskip\labelsep\bf Proposition~\themsct]\it\ignorespaces}{\endtrivlist\vskip3pt}

\newenvironment{m-corollary}{\vskip5pt\refstepcounter{msct}\trivlist \itemindent 0pt%
\item[\hskip\labelsep\bf Corollary~\themsct]\it\ignorespaces}{\endtrivlist\vskip3pt}

\newenvironment{m-lemma}{\vskip5pt\refstepcounter{msct}\trivlist \itemindent 0pt%
\item[\hskip\labelsep\bf Lemma~\themsct]\it\ignorespaces}{\endtrivlist\vskip3pt}

\newenvironment{m-definition}{\vskip5pt\refstepcounter{msct}\trivlist \itemindent0pt%
\item[\hskip\labelsep\bf Definition~\themsct]\ignorespaces}{\endtrivlist\vskip5pt}

\newenvironment{m-notation}{\vskip5pt\refstepcounter{msct}\trivlist \itemindent0pt%
\item[\hskip\labelsep\bf Notation~\themsct]\ignorespaces}{\endtrivlist\vskip5pt}

\newenvironment{m-example}{\vskip5pt\refstepcounter{msct}\trivlist \itemindent0pt%
\item[\hskip\labelsep\bf Example~\themsct]\ignorespaces}{\endtrivlist\vskip5pt}

\newenvironment{m-remark}{\vskip5pt\refstepcounter{msct}\trivlist \itemindent0pt%
\item[\hskip\labelsep\bf Remark~\themsct]\ignorespaces}{\endtrivlist\vskip5pt}

\newenvironment{m-procedure}{\vskip5pt\refstepcounter{msct}\trivlist\itemindent0pt%
\item[\hskip\labelsep\bf Procedure~\themsct]\ignorespaces}{\hfill\endtrivlist\vskip5pt}%

\newenvironment{m-question}{\vskip5pt\refstepcounter{msct}\trivlist \itemindent0pt%
\item[\hskip\labelsep\bf Question.]\ignorespaces}{\endtrivlist\vskip5pt}

\newenvironment{thm-nono}[1]{\vskip5pt\trivlist \itemindent 0pt %
\item[\hskip\labelsep\bf Theorem~{\rm\mbox{#1}}]\it\ignorespaces}{\endtrivlist\vskip5pt}

\newenvironment{lm-nono}[1]{\vskip5pt\trivlist \itemindent0pt%
\item[\hskip\labelsep\bf Lemma~{\rm\mbox{#1}}]\it\ignorespaces}{\endtrivlist\vskip5pt}

\newenvironment{conj-nono}[1]{\vskip5pt\trivlist \itemindent0pt%
\item[\hskip\labelsep\bf Conjecture~{\rm\mbox{#1}}]\it\ignorespaces}{\endtrivlist\vskip5pt}

\newenvironment{m-thank}{\vskip5pt\trivlist \itemindent0pt%
\item[\hskip\labelsep\it Acknowledgments]\ignorespaces}{\endtrivlist\vskip5pt}

\newenvironment{m-proof}{\vskip2pt\trivlist \itemindent0pt%
\item[\hskip\labelsep\it Proof.]\ignorespaces}{\hfill$\Box$\endtrivlist\vskip5pt}%

\newenvironment{m-asmp}{\vskip5pt\trivlist \itemindent0pt%
\item[\hskip\labelsep\bf Assumption.]\ignorespaces}{\hfill\endtrivlist\vskip5pt}%

\newcounter{meqn}[section]\renewcommand{\themeqn}{\thesection.\arabic{meqn}}
\newenvironment{m-eqn}[1]{\vskip5pt\refstepcounter{meqn}%
\trivlist\itemindent0pt\item[]\ignorespaces%
\hfill $\displaystyle #1$\hfill\hbox{\rm(\themeqn)}}{\endtrivlist\vskip5pt}

\newcommand{\bibauth}[2]{\textrm{{#1}~{#2}}}
\newcommand{\bibtitl}[1]{\textrm{#1}.}
\newcommand{\bibjnyp}[4]{\textrm{#1} \textbf{#2} (#3), {#4}.}
\newcommand{\bibinbook}[4]{In: \textrm{#1}\textrm{, #2}\textrm{, #3}\textrm{, #4}.}
\newcommand{\bibbook}[4]{\textrm{#1}. {#2}, {#3}, {(#4)}.}

\let\mt\mapsto

\font\tenmsa=msam10 %
\newcommand\hdashpiece{%
{\vrule height2.75pt depth-2.35pt width2.3pt \kern1.7pt}}%
\newcommand\hdashpieces{%
{\hdashpiece\hdashpiece\hdashpiece\hdashpiece}}%

\newcommand\dashar{\mathrel{%
\hdashpieces\kern-0.4pt\hbox{\tenmsa K}}}%

\DeclareFontFamily{OT1}{rsfs}{}
\DeclareFontShape{OT1}{rsfs}{n}{it}{<->rsfs10}{}
\DeclareMathAlphabet{\crl}{OT1}{rsfs}{n}{it}

\let\disp\displaystyle
\let\ges\geqslant

\let\dta\delta 
\let\les\leqslant
\let\mt\mapsto
\let\nit\noindent

\newcommand\uset[2]{\mathop{#2}\limits_{#1}}

\newcommand{\rd}{{\rm d}}
\newcommand{\ala}{\alpha}
\newcommand{\bta}{\beta}
\newcommand{\bb}{b}
\newcommand{\BB}{B}

\newcommand{\tBB}{{\tilde\BB}}
\newcommand{\tBp}{{\tilde\BB_+}}

\newcommand{\CC}{C}
\newcommand{\cc}{c}
\newcommand{\ee}{\veps}
\newcommand{\het}{{\rm heat}}
\newcommand{\lot}{{l.o.t.}}
\newcommand{\Lda}{{\Lambda}}
\newcommand{\LL}{L}
\newcommand{\resd}{{\rm resd}}

\newcommand{\tup}{{\tilde u}_+}
\newcommand{\tu}{t}
\newcommand{\typ}{\tilde y_+}
\newcommand{\TT}{T}
\let\veps\varepsilon

\title{A study of the One-Dimensional\\ Heat-Conduction Equation with Radiation}
\author{Mihai Halic}


\begin{document}

\begin{abstract}
We consider a boundary value problem (BVP) modelling one-dimensional heat-conduction with radiation, which is derived from the Stefan-Boltzmann law. The problem strongly depends on the parameters, making difficult to estimate the solution. In here we apply an analytical approach to determine upper and lower bounds to the exact solution of the BVP, which allows estimating the latter. Finally, we support our theoretical arguments with numerical data, by implementing them into the MAPLE computer program.
\end{abstract}

\maketitle
\markboth{\sc Mihai Halic}{\sc On heat-conduction with radiation}

\section{Introduction}

Firnett-Troesch~\cite{fi-tr} applied the so-called `shooting method' to numerically analyse two second order BVPs, which depend in a sensitive way on the parameters. The author devised~\cite{hrt} an analytic approach for determining explicit approximate solutions to such BVPs, and compared the outcome with existing numerical results for the `first Troesch-equation' $y''=L\sinh(Ly)$. The analytical method turned out to be several orders of magnitude preciser than previous numerical methods, especially for large values of $L$. 
Here we develop this approach, and consider the `second Troesch-equation' that is, 
$$\hspace{.1\textwidth}
u''_\het=\bb^2(u^4_\het-\tu^4),\quad u_\het(0)=1,\;u_\het(1)=\tu\in(0,1),
\hspace{.1\textwidth}\text{(HCR)}
$$
which describes the steady-state, one-dimensional heat conduction with radiation. The equation depends on two parameters: $b$ is related to the radiative properties of materials, while $t$ is roughly the ratio between the initial and final temperatures. 

Our goal here is to approximately solve (HCR), understand the dependence of its solution on the parameters, elucidate the law governing the boundary layer, etc. These objectives turn out to be challenging. There seems to be a rather restricted literature dedicated to the analysis of this very  equation, in spite of being interesting:
\begin{itemize}
\item On the mathematics side, the ODE is far from being classically integrable. Existence results are proved  in~\cite{hu-ca,ma-wo,qly}. Power series expansions are obtained in~\cite{tao,tok}; their terms are complicated and, for good approximations, the order of the development quickly increases with the parameters. 

From a numerical perspective, (HCR) often has the so-called `boundary layer property', $u_\het$ has huge derivative near $x=0$, which is a major obstacle. 
For running the shooting method, one needs a very good estimate of either $u'_\het(0)\ll0$ or $u'_\het(1)\approx0$. Unfortunately, for larger $\bb$ ---physically realistic situations involve $\bb\approx30 \dots 40$--- these values are hard to estimate. Indeed, Firnett-Troesch consider only a few values ($\,0.8, 0.4, 0\,$) of $\tu$, and deal with values of $\bb$ up to $20$. \smallskip 

\begin{center}\begin{tabular}{cc}
\scriptsize \cite[p.427]{fi-tr}:&\scriptsize \cite[Figure 6.1, p.426]{fi-tr}:\\ 
\begin{minipage}{.355\textwidth}
\scriptsize ``If the differential equation is integrated with a guess [...] which is even slightly in error, then a strongly diverging solution is obtained." 
\end{minipage}
&
\begin{minipage}{.355\textwidth}
\centering{
\includegraphics[width=.95\textwidth]{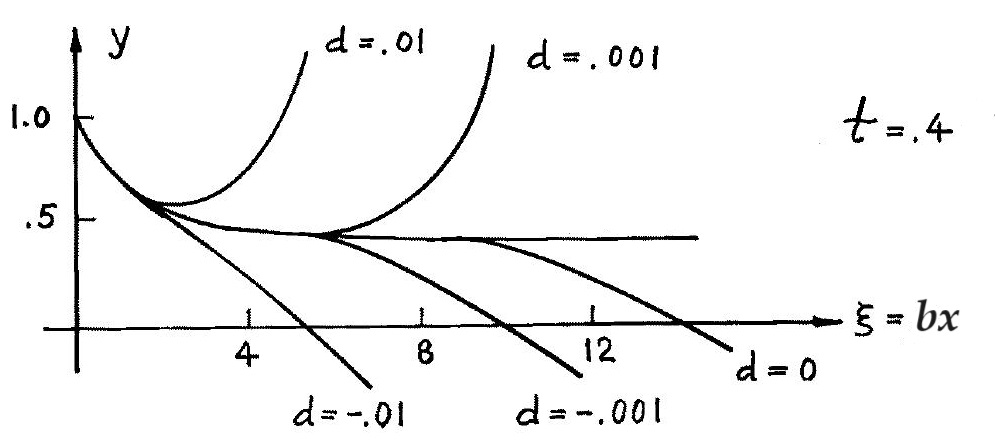}}
\end{minipage}
\end{tabular}\end{center}
\smallskip 

Thus one faces a chaotic situation ---tiny changes lead to major variations---, which raises the issue how to deal with a non-integrable ODE, without completely relying on computer power. 
Our approach is to determine explicit upper and lower bounds ---so-called envelopes--- and deduce the dependence on the parameters. Quite remarkably, the conclusions follow naturally.

\item On the physics side, the equation is based on the Stefan-Boltzmann law for the heat radiation; it's a rare physical law involving $4^{th}$-power dependence on a parameter (the temperature). Numerous textbooks~\cite{arp,ca-ja,can,hol,mod}, written from both theoretical and engineering standpoints, show the importance of heat conduction and radiation. 
\end{itemize}


\subsection*{Outline}

At the core of our approach, there is a careful analysis of the heat  conduction equation. The large variations of the output, generated by insignificant changes in the input, can't be controlled without taming the mechanism which yields solution. The computer work supplements and enhances the presentation, also reinforces the (possibly confusing) theoretical arguments with numerical evidence. 
 
The article is organized as follows: we start by estimating the solution of a transcendental equation, which occurs in our analysis. In Section~\ref{sct:heat1} we obtain the first estimates for $u_\het$, actually an upper bound. This estimate is already {\em crucial}: 
\begin{itemize}
\item It yields an approximate value for the initial derivative, valid for {\em any choices} of the parameters $\bb$ and $\tu$, which is otherwise impossible to guess. This approximate value allows running the Runge-Kutta method for (HCR), with quite a remarkable precision, even for low values of the parameters.\footnote{The author reported on the application of the Runge-Kutta method, briefly discussed here in Section~\ref{ssct:shoot}, at the ICNAAM~2023 conference. The short abstract of the report will appear in the AIP Conference Proceedings.}  
\item It allows pinning down the law governing the boundary layer property, in terms of the parameters $\bb,\tu$. Remarkably, this law has a geometric and a variational interpretation.
\end{itemize}
In the Sections~\ref{sct:heat2},~\ref{sct:bounds},~\ref{sct:back-heat}, we refine our considerations, and determine accurate envelopes for the exact solution of the heat equation. 
Finally, we implement our considerations into the MAPLE software, in Section~\ref{sct:num}; the code is given, the interested reader may freely experiment with the favourite values of the parameters. Numerical data shows that the precision of the estimates increases with the slope of the boundary layer.

\begin{m-thank}
The author thanks K. Mustapha for discussions related to numerical modelling, which helped improving the rule governing the boundary layer.
\end{m-thank}


\section{Background matter}\label{sct:back}

The following calculations are tedious, but essential. 

\subsection{A transcendental equation\;I}\label{ssct:treq-1}

We consider real numbers $\zeta>1$, $\Lda>0$, and the equation in the unknown $\cc\ges0$: 
\begin{m-eqn}{
\begin{array}{r}
\disp\frac{\zeta-\cc}{\zeta+\cc}=\frac{1-\cc}{1+\cc}\cdot e^{2\Lda \cc}
\;\Leftrightarrow\;
(\zeta-1)\cdot\frac{\cc}{\zeta-\cc^2}=\tanh(\Lda \cc)
\;\Leftrightarrow\;
\zeta= \underbrace{\cc\cdot\frac{1-\cc\cdot \tanh(\Lda \cc)}{\cc-\tanh(\Lda \cc)}}_{rhs(\cc)}.
\end{array}
}\label{eq:cc}
\end{m-eqn}
First we analyse the poles of the right-hand side. The equation $\cc=\tanh(\Lda \cc)$ has the solution $\cc_0(\Lda)=0$, for $\Lda\in(0,1]$. For $\Lda>1$, it admits a (unique) solution $\cc_0(\Lda)\in(0,1)$. In all cases, $\cc_0(\Lda)$ can be expressed as an infinite iteration of the `tanh' function: 
$$
\cc_0(\Lda)=\tanh(\Lda\cdot\tanh(\Lda\cdot\tanh(\Lda\cdot \dots))).
$$
This formula allows computing approximate values for $\cc_0(\Lda)$:
\begin{m-eqn}{
\begin{array}{|l|l|}
\hline 
\Lda&\cc_0(\Lda)=\tanh(\Lda\cdot\tanh(\Lda\cdot\tanh(\Lda\cdot \dots))) \\ 
\hline 
\Lda=1& \cc_0(1)=0 \\ 
\hline 
\Lda=1.10& \cc_0(1.10)\approx 0.5\\ 
\hline 
\Lda=1.50& \cc_0(1.50)\approx 0.85\\ 
\hline 
\Lda=1.67& \cc_0(1.67)\approx 0.9\\ 
\hline 
\Lda=2& \cc_0(2)\approx 0.95\\ 
\hline 
\Lda\ges5& \cc_0(\Lda)>0.999\approx 1\\ 
\hline 
\end{array}
}\label{tab:cc}
\end{m-eqn}
A numerical experiment shows that $20$ iterations of the function $\tanh(\Lda\cdot)$ are necessary to obtain a two digit approximation ---an upper bound--- of $\cc_0(\Lda)$, for $\Lda\ges1.1$. 

\begin{m-lemma}\label{lm:cc}
Given $\zeta>1$ and $\Lda>0$, the equation~\eqref{eq:cc} admits a unique solution $c=c(\zeta,\Lda)$, as follows: 
\begin{enumerate}
\item 
for $\Lda>1$,  one has $\cc\in(\cc_0(\Lda),1)$.
\item 
for $\Lda\in(0,1]$, $\zeta<\frac{1}{1-\Lda}$, one has $\cc\in(0,1)$.  
(For $\Lda=1$, we let $1/(1-1)=+\infty$.)
\end{enumerate}
\end{m-lemma}

\begin{proof}
In the first case, we have $\uset{\cc\to \cc_0(\Lda)^+}{\lim}rhs(\cc) =+\infty$ and $rhs(1)=1<\zeta$. In the second case, the denominator doesn't vanish, and we have $\uset{\cc\to 0^+}{\lim}rhs(\cc) =\frac{1}{1-\Lda}>\zeta$, $rhs(1)=1<\zeta$. Since $c\mt rhs(c)$ is decreasing, the uniqueness follows. 
\end{proof}


\subsection{Bounding the solution \textit{c}}\label{ssct:bound-c}

We determine explicit upper and lower bounds for the exact solution of~\eqref{eq:cc}, depending on the parameters $\zeta,\Lda$. The previous Lemma shows that we should distinguish between two cases: $\Lda>1$ (actually, we take $\Lda\ges1.1$) and $\Lda<1.1$. The approaches to estimating the solutions are different.


\subsubsection{Case $\Lda\ges1.5$} 

Polynomial approximations fail in this situation because the $\tanh$-function approaches exponentially fast the value $1$. Since $\cc$ belongs to $(0,1)$, we write: 
$$\;\cc=\tanh(Y)=\tanh(\Lda X)>\cc_0(\Lda)=\tanh(\Lda\cdot \cc_0(\Lda)),\; X>\cc_0(\Lda).$$ 

The addition formula for the $\tanh$-function yields:
\begin{m-eqn}{
\zeta=\frac{\cc}{\tanh(Y-\Lda \tanh(Y))}=\frac{\tanh(\Lda X)}{\tanh(\Lda(X-\tanh(\Lda X)))},\quad X=\cc_0(\Lda)+\gamma,\;\gamma>0.
}\label{eq:X}
\end{m-eqn}
To determine explicit approximations, we consider the inequalities `$\les$' and `$\ges$'. 

\begin{itemize}
\item[`$\les$'] 
Since $\cc>\cc_0(\Lda)$, it's enough to have $\tanh\big(\Lda (X-\tanh(\Lda X))\big)\les\zeta^{-1}\cc_0(\Lda)$, so 
$$
X-\tanh(\Lda X)\les \Lda^{-1}\cdot\tanh^{-1}(\zeta^{-1}\cdot \cc_0(\Lda)).
$$
By using~\eqref{eq:X}, the left-hand side is less than $\gamma+\cc_0(\Lda)-\tanh(\Lda \cc_0(\Lda))=\gamma$. Overall, we deduce that the following value is convenient:
\begin{m-eqn}{
\cc_-(\zeta,\Lda)=\tanh\biggl( \Lda\cdot \cc_0(\Lda) + \tanh^{-1}(\zeta^{-1}\cdot \cc_0(\Lda)) \biggr). 
}\label{eq:ccm}
\end{m-eqn}
\item[`$\ges$'] 
As $c<1$, it's enough to have: 
$$
\begin{array}{l}
\tanh\big(\Lda (X-\tanh(\Lda X))\big)\ges\zeta^{-1}\ges\zeta^{-1}\cc,
\\[2ex] 
\gamma+\tanh(\Lda \cc_0(\Lda))-\tanh(\Lda X)=X-\tanh(\Lda X)\ges\Lda^{-1}\tanh^{-1}(\zeta^{-1}).
\end{array}
$$
A concavity argument yields $\tanh(\Lda X)-\tanh(\Lda \cc_0(\Lda))\les\frac{\Lda\gamma}{\cosh^2(\Lda \cc_0(\Lda))}$, which implies that the following value is convenient:
\begin{m-eqn}{
\cc_+(\zeta,\Lda)= \tanh\biggl(\Lda\cdot \cc_0(\Lda) + \frac{\tanh^{-1}(\zeta^{-1})}{1-\frac{\Lda}{\cosh^2(\Lda\cdot \cc_0(\Lda))}}\biggr).
}\label{eq:ccp}
\end{m-eqn}
Note that, by differentiating $\cc_0(\Lda)=\tanh(\Lda \cc_0(\Lda))$ for $\Lda$, one finds that the denominator above is positive.
\end{itemize}


\subsubsection{Case $\Lda\ges5,\;\zeta\les e^{\Lda}$}

In this situation, the previous formulae can be simplified, by eliminating $c_0(\Lda)\approx1$: we replace $X=1+\gamma$ in~\eqref{eq:X}, with $\gamma\ges0$. The very same approach yields: 
\begin{m-eqn}{
\begin{array}{rl}
\text{`$\les$' is satisfied by}&\disp 
c_-(\zeta,\Lda):=\tanh\Big(\Lda\cdot\tanh(\Lda)+\tanh^{-1}\big(\zeta^{-1}\cdot\tanh(\Lda) \big) \Big);
\\[2ex] 
\text{`$\ges$' is satisfied by}&\disp 
c_+(\zeta,\Lda):=\tanh\Biggl(
\frac{\tanh^{-1}(\zeta^{-1})+\frac{\Lda(\sinh(2\Lda)-2\Lda)}{2\cosh^2(\Lda)}}
{1-\frac{\Lda}{\cosh^2(\Lda)}}  \Biggr).
\end{array}
}\label{eq:cc5}
\end{m-eqn}
The condition $\zeta\les e^\Lda$ is necessary to ensure that $\gamma$ is positive. 


\subsubsection{Case $\Lda\in(0,1.5)$} 

In this case, the bounds $c_0(\Lda)< c< 1$ are imprecise, we need better ones. The series development of $\frac{\tanh(x)}{x}$ yields the estimate: 
\begin{m-eqn}{
\begin{array}{c}
\Bigl|
\frac{\tanh(x)}{x}-(1-0.2x^2)
\Bigr|\les 0.15x^2,\quad\forall\,x\in[0,1.5].
\end{array}
}\label{eq:tanh-bound}
\end{m-eqn}
Thus the solutions of the bi-quadratic equations 
\begin{m-eqn}{
\begin{array}{l}
\frac{\zeta-1}{\Lda}=(\zeta-\cc^2)\cdot\Bigl(1 - 0.35\cdot\Lda^2\cc^2 \Bigr),\quad 
\frac{\zeta-1}{\Lda}=(\zeta-\cc^2)\cdot\Bigl(1 - 0.05\cdot\Lda^2\cc^2 \Bigr),
\\[2ex]\disp
\cc_\rho(\zeta,\Lda)=
{\Biggl[ \frac{2\frac{1-\zeta(\Lda-1)}{\Lda}}{(1+\rho\cdot\zeta\Lda^2)+\sqrt{(1-\rho\cdot\zeta\Lda^2)^2+4\rho\Lda(\zeta-1)}}
\Biggr]}^{1/2},\; \rho=0.35\;and\;0.05,
\end{array}
}\label{eq:c-Lsmall}
\end{m-eqn}
are bounds $\cc_{0.35}(\zeta,\Lda)<\cc<\cc_{0.05}(\zeta,\Lda)$ to the exact solution of~\eqref{eq:cc}. We use these values in~\eqref{eq:X}; for shorthand, we drop the variables $(\zeta,\Lda)$. 
\begin{itemize}
\item[`$\les$'] It's enough to have  $\tanh\big(\Lda (X-\tanh(\Lda X))\big)\approx\zeta^{-1}\cc_{0.35}$, so we set: 
$$
c_-(\zeta,\Lda):=\tanh\Big(\Lda\cdot\cc_{0.35}+\tanh^{-1}\big(\zeta^{-1}\cdot\cc_{0.35}\big)\Big).
$$
\item[`$\ges$'] It's enough to have $\tanh\big(\Lda (X-\tanh(\Lda X))\big)\approx\zeta^{-1}\cc_{0.05}$, so we set:  
$$
c_+(\zeta,\Lda):=
\tanh\Biggl(
\frac{\tanh^{-1}(\zeta^{-1}\cc_{0.05})+\frac{\Lda\big(\sinh(2\Lda\cc_{0.05})-2\Lda\cc_{0.05}\big)}{2\cosh^2(\Lda\cc_{0.05})}}
{1-\frac{\Lda\cc_{0.05}}{\cosh^2(\Lda\cc_{0.05})}}  \Biggr).
$$
\end{itemize}


\begin{m-remark}\label{rmk:cc} 
(i) One should keep in mind that $\cc_+(\zeta,\Lda)\les c\les c_-(\zeta,\Lda)$. The notation is chosen in such a way that the $-$ (resp. $+$) indices determines lower (resp. upper) bounds for the exact solution of the heat equation. 

\nit(ii) The equations~\eqref{eq:ccm},~\eqref{eq:ccp} involve $\cc_0(\Lda)$, which is implicitly defined. For numerical computations, we replace it by $30$ iterations of $\tanh(\Lda\cdot)$ that is:
$$
c_0(\Lda)\approx\underbrace{\tanh(\Lda\cdot\tanh(\dots(\tanh(\Lda)\dots)))}_{30\;\text{times}},\quad\Lda\in[1.5,5].
$$
For large values of $\Lda$, one may use the estimates~\eqref{eq:cc5}. 
\label{rmk:treq-2}
\end{m-remark}


\subsection{A transcendental equation\;II} \label{ssct:treq-2}

For real numbers $z_1>z_0>0$, $\LL>0$, $\bta\ges\ala>0$, we consider the equation in $\CC\ges0$, 
\begin{m-eqn}{
\frac{z_1-\ala\CC}{z_1+\bta\CC}=\frac{z_0-\ala\CC}{z_0+\bta\CC}\cdot e^{2\LL\CC}. 
}\label{eq:zetac}
\end{m-eqn}
(For our purposes, we'll need $\bta=z_0=1$, only.) The substitutions 
$$
\zeta_i:=z_i+\frac{\bta-\ala}{2}\CC\;(i=0,1),\quad\cc:=\frac{\bta+\ala}{2\zeta_0}\CC,\quad\Lda:=\frac{2\zeta_0\LL}{\bta+\ala}, \quad \zeta:=\frac{\zeta_1}{\zeta_0},
$$ 
reduce this equation to~\eqref{eq:cc}, but extra complications occur. The unknown $\CC$ is given by the equation
\begin{m-eqn}{
\CC=\frac{2z_0\cc}{(\bta+\ala)-(\bta-\ala)\cc},
}\label{eq:CC}
\end{m-eqn}
so $\zeta,\cc$ are related.  (Before, $\cc$ and $\zeta$ were independent.) Fortunately, we are interested in solving the inequalities `$\les$' and `$\ges$' in~\eqref{eq:zetac}, which lead to 
$$
\zeta=\frac{z_1+\frac{\bta-\ala}{2}\CC}{z_0+\frac{\bta-\ala}{2}\CC}
=\cc\cdot\frac{1-\cc\cdot\tanh(\Lda\cc)}{\cc-\tanh(\Lda\cc)}. 
$$
with inequalities in the same direction, respectively. 
\begin{itemize}
\item The left-hand side is decreasing in the variable $\CC$ (see~\eqref{eq:CC}), and we have 
\begin{m-eqn}{
\frac{z_0c_0(\Lda)}{\beta}\les\frac{2z_0\cc_0(\Lda)}{(\bta+\ala)-(\bta-\ala)\cc_0(\Lda)}<\CC<\frac{z_0}{\ala}.
}\label{eq:CC-bounds}
\end{m-eqn}
\item The right-hand side is  increasing in $\Lda$, and we have $\frac{z_0}{\bta}L\les\Lda\les\frac{z_0}{\ala}L$.
\end{itemize}
Thus we obtain approximate solutions of the inequalities:
\begin{itemize}
\item[`$\les$'] It requires solving  
$$
\zeta_-:=\frac{z_1+\frac{(\bta-\ala)z_0}{2\bta}c_0(\frac{z_0}{\beta}L)}{z_0+\frac{(\bta-\ala)z_0}{2\bta}c_0(\frac{z_0}{\beta}L)}\les 
\cc\cdot\frac{1-\cc\cdot\tanh(\frac{z_0}{\bta}L\cc)}{\cc-\tanh(\frac{z_0}{\bta}L\cc)},
$$
which is satisfied by $\cc=\cc_-(\zeta_-,\frac{z_0}{\bta}L)$. 
\item[`$\ges$'] It requires solving 
$$
\zeta_+:=\frac{z_1+\frac{(\bta-\ala)z_0}{2\ala}}{z_0+\frac{(\bta-\ala)z_0}{2\ala}}\ges
\cc\cdot\frac{1-\cc\cdot\tanh(\frac{z_0}{\ala}L\cc)}{\cc-\tanh(\frac{z_0}{\ala}L\cc)},
$$ 
which is satisfied by $\cc=\cc_+(\zeta_+,\frac{z_0}{\ala}L)$. 
\end{itemize}
The corresponding values $C_\pm$ are obtained by inserting $c_\pm$ into~\eqref{eq:CC}.


\subsection{Solution to an ODE}\label{ssct:ode}

Some calculations yield the following:

\begin{m-proposition}\label{prop:ode}
Given $0<\veps<1$ and $0<C<\veps^{-1}$. The solution to the ODE 
$$
\frac{y'}{\sqrt{2}}=\frac{2\tBB}{1+\ee}\bigl(y^{2.5}+(1-\ee)\CC y-\ee\CC^2y^{-0.5}\bigr),\quad y(1)=\TT, 
$$
is the function 
\begin{m-eqn}{
y(x)=y_{\tBB,\CC,\ee}(x):=
\biggl[ 
\CC\cdot\frac{\TT^{1.5}\cdot (1+\ee\cdot e^{3\sqrt{2}\tBB\CC(1-x)})+\ee \CC(e^{3\sqrt{2}\tBB\CC(1-x)}-1)}
{T^{1.5}(e^{3\sqrt{2}\tBB\CC(1-x)}-1)+\CC(e^{3\sqrt{2}\tBB\CC(1-x)}+\ee)}
 \biggr]^{2/3}.
}\label{eq:ODE}
\end{m-eqn}
This function has the boundary value $y(0)=1$ precisely when $\tBB,\CC,\ee$ satisfy the equation~\eqref{eq:zetac}, with $z_0=1,z_1=T,\alpha=\veps,\beta=1$.
\end{m-proposition}

\subsection{Criteria for upper/lower envelopes} 

We will determine global, as well as partial envelopes for $u_\het,y_\het$. For deciding that a function is indeed a global envelope, we'll use the following criterion:

\begin{m-proposition}\label{prop:+-}
Let $f,g$ be continuous, positive on $[m,M]$, such that one of them is strictly increasing. We consider the solutions to the following boundary-value problems: 
$$
\begin{array}{llr}
&y_0''=f(y_0),&y_0(0),y_0(1)\;given\;in\;[m,M],\\
&y''=g(y),&y(0),y(1)\;given\;in\;[m,M].\\
\end{array}
$$ 
If one has the following three inequalities 
$$g-f,\; y_0(0)-y(0),\;y_0(1)-y(1)\ges0,\quad\text{(resp. $\les0$)}$$ 
then we have $y_0-y\ges0$ (resp. $\les0$).
\end{m-proposition}

In this article, we'll consider $f(y)=B^2(5y^4-5)$. The proposition can be reformulated in terms of the  so-called `residue functional' corresponding to it, which associates the function  
$$
\resd(y):=y''-B^2(5y^4-5), 
$$
to any (differentiable) function $y$. The residue measures the deviation of $y$ from being a solution to the given ODE. Our criterion reads as follows: 
\begin{itemize}
\item If $y(0)\ges1, y(1)\ges T$, and $\resd(y)\les0$, then $y$ is greater than $y_\het$.
\item If $y(0)\les1, y(1)\les T$, and $\resd(y)\ges0$, then $y$ is smaller than $y_\het$.
\end{itemize}

This is a criterion for global envelopes. However, we'll be also interested in partial envelopes, along the interesting boundary layer. 

\begin{m-proposition}\label{prop:part+-}
Let $F,G$ be positive, differentiable on $[1,T]$, and consider the IVPs 
$$
\begin{array}{llll}
y_0'=\sqrt{2}\cdot F(y_0),&\;y'=\sqrt{2}\cdot G(y),
&&y_0(1)=y_0(1)=T\;\text{given}.
\end{array}
$$
If the following inequality is satisfied (where $r\in[1,T]$):
$$
G-F\ges0,\; \text{for}\;y\in[r,T],\quad \text{(resp. $\les0$)},
$$
then $y_0- y\ges0$ (resp. $\les0$), whenever $y\in[r,T]$. 
\end{m-proposition}

The function of interest will be  $F(y)=B^2\cdot\sqrt{y^5-5y+4+\dta}$, with $\dta\ges0$. Let us explain how this will be used to determine partial envelopes in Propositions~\ref{prop:partial-} and~\ref{prop:partial+}. 
The IVPs are obtained by respectively integrating the BVPs, so 
$$
F(y)=\sqrt{\dta+\int_0^yf(\xi)\rd\xi},\;\text{ and similarly for}\;G.
$$
Suppose we are in the following situation: $G(r)\ges F(r)$ and there are solutions to the IVPs above, such that $\resd(y)=g(y)-f(y)\ges0,\;y\in[r,T]$; this condition holds, in particular, if  $\resd(y=r)\ges0$ and $\resd(y)$ in increasing on $[r,T]$).  Then Proposition~\ref{prop:part+-} applies.


\section{The heat conduction equation: I}\label{sct:heat1}

\subsection{Estimates} 

We are interested in bounding the solution of the BVP:
\begin{m-eqn}{
u''_\het=\bb^2(u^4_\het-\tu^4),\quad u_\het(0)=1,\;u_\het(1)=\tu\in(0,1).
}\label{eq:uhet}
\end{m-eqn}
We perform several changes of variables, to bring the problem into a more suitable form. The substitution $x\mt 1-x$ yields 
\begin{m-eqn}{
u''=\bb^2(u^4-\tu^4),\quad u(0)=\tu\in(0,1),\;u(1)=1.
}\label{eq:ubt}
\end{m-eqn}
By integrating, we obtain 
$$
\frac{(u')^2}{2\tu^2}=\frac{\bb^2\tu^3}{5}\Bigl(\frac{u^5}{\tu^5}-5\frac{u}{\tu}+A\Bigr).
$$
The derivative at the left-end of the interval $[0,1]$ is expected to be close to zero. We set $A=4+\dta$ and $y:=u/\tu$, so the equation becomes:
\begin{m-eqn}{
\begin{array}{ll}\disp 
\frac{y'_\het}{\sqrt{2}}=\BB\cdot\sqrt{y^5_\het-5y_\het+4+\dta},
&\disp y_\het(0)=1,y_\het(1)=\tu^{-1},
\\[2ex]\disp 
y''_\het=5\BB^2(y^4_\het-1),&\disp \BB:=\frac{\bb\tu^{3/2}}{\sqrt{5}}.
\end{array}
}\label{eq:y-het}
\end{m-eqn}
Clearly, the equation can't be integrated exactly, so we need to approximate it. At the left-end of $[0,1]$, we have $y_\het(0)=1$. However, at the right-end, the function $y_\het$ take large values, for $\tu$ close to zero, so we need to take into account the large-$y$ behaviour. Therefore, we consider the boundary value problem
$${
\frac{y'}{\sqrt{2}}=\BB\cdot G(y),\quad y(0)=1,y(1)=\frac{1}{\tu},
}$$
with the right-hand side satisfying the following requirements: 
\begin{itemize}
\item The difference $\sqrt{y^5-5y+4+\dta}-G(y)$ is small, for all $y\ges1$;
\item The equation $y'=\BB\cdot G(y)$ is explicitly integrable. 
\end{itemize}

For large $y$, we have the expansion 
$$
\sqrt{y^5-5y+4+\dta}
=y^{\frac{5}{2}}-\frac{5}{2}y^{-\frac{3}{2}}+\frac{4+\dta}{2}y^{-\frac{5}{2}}+\dots
=y^{\frac{5}{2}}\cdot\biggl( 1- \frac{5}{2}y^{-4} +\frac{4+\dta}{2}y^{-5}+\cdots \biggr), 
$$
which hints to considering $G(y)=y^{5/2}\cdot c(y)$, such that the approximation is still good for $y\approx1$. Thus the function $c(y)$ must satisfy $c(1)\approx 0$ and $\uset{y\to\infty}{\lim} c(y)=1$.  The asymptotic expansion hints to considering  
$$
G(y)=y^{2.5}\cdot(1-\text{const.}\cdot y^{-2.5-r})=y^{2.5}-\text{const.}\cdot y^{-r},
$$
with $\text{const.}\approx 1$. 
It remains to determine a value of $r$ which is appropriate for both requirements. The numerator of the difference  
$$
\sqrt{y^5-5y+4} - (y^{5/2}-y^{-r})=\frac{-5y+4+2y^{2.5-r}-y^{-2r}}{\sqrt{y^5-5y+4}+(y^{5/2}-y^{-r})}, 
$$
vanishes of second order at $y=1$ for $r_0=\sqrt{10}-2.5\approx0.6623$; for approximation purposes, a convenient value for $r$ should be chosen nearby. However, the resulting ODE must be explicitly integrable, too, and this determines $r=0.5$.


\subsection{A boundary value problem\,I}\label{ssct:bvp1} 

The value $r=0.5$ leads to the following BVP (where we set $T:=\tu^{-1}$): 
\begin{m-eqn}{
\begin{array}{l|l}
\frac{y'}{\sqrt{2}}=\BB\cdot(y^{2.5}-\CC^2y^{-0.5}),&y''=\BB^2(5y^4-4\CC^2y-\CC^4y^{-2}),
\\ 
y(0)=1, y(1)=T,&y(0)=1, y(1)=T.
\end{array}
}\label{eq:heat}
\end{m-eqn}
In the first formulation, the constant $\CC>0$ is determined by the boundary conditions. For approximating the exact solution, we change our point of view. We drop the boundary value at $x=0$, and consider the initial value problem (IVP):
\begin{m-eqn}{
\begin{array}{l}
\frac{y'}{\sqrt{2}}=\tBB\cdot(y^{2.5}-\CC^2\cdot y^{-0.5}),\quad y\ges1,\; \CC\in(0,1],
\\[1ex]
y(1)=T,
\end{array}
}\label{eq:odey}
\end{m-eqn}
whose general solution is: 
\begin{m-eqn}{
\begin{array}{l}
\disp
\frac{T^{1.5}-\CC}{T^{1.5}+\CC}=\frac{y^{1.5}-\CC}{y^{1.5}+\CC}\cdot e^{3\sqrt{2}\tBB\CC\cdot (1-x)},\quad 
\frac{T^{1.5}-\CC}{T^{1.5}+\CC}
=\frac{1-\CC}{1+\CC}\cdot e^{3\sqrt{2}\tBB\CC},
\\[2ex]
\disp 
y={\Biggl[ 
\frac{C}{\tanh\Big( \tanh^{-1}(\CC t^{1.5})+1.5\sqrt{2}\tBB\CC(1-x) \Big)} 
\Biggr]}^{2/3}\kern-1ex.
\end{array}
}\label{eq:y}
\end{m-eqn}
(The parameter $\tBB$ should be thought off as an `improved value' for $\BB$, which increases the accuracy of the estimates.) The equation defining $\CC$ was analysed in Section~\ref{sct:back}, with the parameter $\Lda=1.5\sqrt{2}\tBB$; by our discussion, we have $\CC\les1$. To obtain an upper bound for~\eqref{eq:y-het}, the second derivative should satisfy Proposition~\ref{prop:+-}:
\begin{m-eqn}{
\begin{array}{l}
y''=\tBB^2(5y^4-4\CC^2y-\CC^4y^{-2})\les\BB^2(5y^4-5),\quad\forall\,y\in[1,\TT],\quad y(0)\ges1.
\end{array}
}\label{eq:tBB}
\end{m-eqn}
The left-hand side is decreasing in $\CC$, so it's enough to ensure that the inequality is satisfied for $\CC=1$. Let 
\begin{m-eqn}{
\tBp:=(\tilde q_+)^{-1}\cdot\BB,\;\text{with}\;
\tilde q_+:=\Bigl[\frac{(1-\tu^{3})(5+\tu^{3})}{5(1-\tu^4)}\Bigr]^{1/2}\in[\sqrt{0.9},1].
}\label{eq:tBp}
\end{m-eqn}

\begin{m-lemma}\label{lm:tBp}
The function 
$\disp \typ(x):= {
\Big[\tanh\Bigl(\tanh^{-1}(\tu^{1.5})+1.5\sqrt{2}\tBp(1-x)\Bigr)\Big]}^{-\frac{2}{3}},
$ 
satisfies 
$$
\begin{array}{l}
\typ''= {(\tBp)}^2(5\typ^4-4\typ-\typ^{-2}),
\\[2ex]
\typ(1)=T=\tu^{-1},\;
\typ(0)= \biggl(1+\frac{2}{\exp\bigl(2\tanh^{-1}(\tu^{1.5})+3\sqrt{2}\tBp\bigr)-1}\biggr)^{\frac{2}{3}},
\end{array}
$$ 
and it is greater than the exact solution $y_\het$ of the BVP~\eqref{eq:y-het}. 
\end{m-lemma}

\begin{m-proof}
The function $\typ$ satisfies the conditions of Proposition~\ref{prop:+-}. 
\end{m-proof}

We briefly pause the analysis to analyse the behaviour of the function $y_\het$, in order to justify the subsequent considerations and refinements.


\subsection{Estimating the derivative}

The major issue for the shooting method~\cite{fi-tr} is the difficulty to estimate the initial derivative $y'_\het(0)=\sqrt{\dta}$. A crucial by-product of Lemma~\ref{lm:tBp} is an upper bound for $\dta$ and, a posteriori, for the derivative at $x=1$.  

\begin{m-proposition}\label{prop:dery1}
The following estimates hold for the derivative of $y_\het$: 
\begin{m-eqn}{
\scalebox{0.85}{$
\begin{array}{l}
\disp\dta\les\Delta:={\Bigg[\frac{\typ(0.25)-1}{0.25\sqrt{2}\BB}\Bigg]}^2
\les \frac{15}{\BB^2}\cdot\frac{1}{(e^{3\BB}-1)^2}\Big(1+\frac{1}{e^{3\BB}-1}\Big)^2;
\\ 
\disp 
0\les y'_\het(1)-\sqrt{2}\BB\sqrt{T^5-5T+4}
\les 
\frac{\sqrt{2}\BB\dta}{2\sqrt{\TT^5-5\TT+4}}
\les 
\frac{11\Big(1+\frac{1}{e^{3\BB}-1}\Big)^2}{\BB(e^{3\BB}-1)^2\cdot\sqrt{\TT^5-5\TT+4}}.
\end{array}
$}
}\label{eq:dery1}
\end{m-eqn}
\end{m-proposition}
The moral is that, already for $\BB\ges2$, the derivative $y_\het'(1)$ differs from the (universal) lower bound $\sqrt{2}\BB\sqrt{\TT^5-5\TT+4}$ by a tiny amount, of size $o(\BB^{-1}e^{-6\BB}\TT^{-2.5})$. Indeed, for  $\BB\ges\frac{1}{3}$, one has $\dta\les 100\BB^{-2}e^{-6\BB}.$

\begin{m-proof}
Note that $y_\het$ is convex ---that is, $y_\het''>0$---, so its derivative at $x=0$ is less than the slope of the secant line between $x=0$ and $x=h=0.25$. Since $y_\het$ is less than $\typ$, we deduce the following estimate: 
$$
\begin{array}{rl}
\sqrt{2}\BB\sqrt{\dta}= y_\het'(0)
&
\les \frac{\typ(h)-1}{h}\les\frac{\typ(h)^3-1}{3h}
=\frac{2\big(\typ(h)^{3/2}-1\big) +\big(\typ(h)^{3/2}-1\big)^2}{3h}.
\end{array}
$$
By inserting $x=h$ into the defining formula, we obtain $\typ(h)^{3/2}-1\les\frac{2}{e^{3.1\BB}-1}$, thus 
$$
\sqrt{\dta}
\les
\frac{4}{h\cdot3\sqrt{2}\cdot\BB}\cdot\frac{1}{e^{3.1\BB}-1}\Big(1+\frac{1}{e^{3\BB}-1}\Big).
$$
It remains to use this estimate, to bound the derivative of $y_\het$ at $x=1$.
\end{m-proof}


\subsection{Application to the `shooting' method}\label{ssct:shoot}

The content of this section should probably belong to Section~\ref{sct:num}, which is devoted to numerical simulations. However, we believe that this short digression is useful to present some applications following from theoretical considerations. 

The shooting method is doomed to fail without an accurate knowledge of the initial derivative of $y_\het$. According to Proposition~\ref{prop:dery1}, we have 
$$
0\les\dta\les\Delta={\Bigg[\frac{\typ(0.25)-1}{0.25\sqrt{2}\BB}\Bigg]}^2.
$$
This bare inequality already yields three possible ways of approximating $y_\het$ (the abbreviation `RK' stands for Runge-Kutta):
\begin{enumerate}
\item Let $\dta=0$, and consider the IVP $y'/\sqrt{2}=\BB\sqrt{y^5-5y+4},\; y(1)=\TT$, whose solution is denoted by $y_{RK+}$. Its second derivative satisfies the same equation $y''=5\BB^2(y^4-1)$ as $y_\het$, and $y_{RK+}'(1)<y_\het'(1)$, hence $y_{RK+}>y_\het$. 
\item Let $\dta=\Delta$, and consider $y'/\sqrt{2}=\BB\sqrt{y^5-5y+4+\Delta},\; y(1)=\TT$, whose solution is denoted by $y_{RK-}$. Its second derivative satisfies $y''=5\BB^2(y^4-1)$, same as $y_\het$, and $y_{RK-}'(1)>y_\het'(1)$, hence $y_{RK-}<y_\het$. 
\item Let $\dta=d:={\bigg[\frac{\typ(0.25)-\typ(0)}{0.25\sqrt{2}\BB}\bigg]}^2$, which is between $0$ and $\Delta$, and consider the corresponding IVP. Its solution $y_{RK}$ will approximate $y_\het$ the best. 
\end{enumerate}
One is led to asking what is the precision of these approximations.

\begin{m-corollary}
\nit The error of the approximation $y_\het(x)\approx  y_{RK}(x)$ is at most $err_y(x)=y_{RK+}(x)-y_{RK-}(x)$. The size of the maximal error equals $O(e^{-3\BB})$.

\nit The error for $u_\het(1-x)\approx \tu\cdot y_{RK}(x)$ has the order $O(t\cdot e^{-3B})$.
\end{m-corollary}

\begin{m-proof}
The error function is decreasing, has negative derivative, so the maximum is attained at $x=0$, and is bounded above by $const\cdot\sqrt{\Delta}$.
\end{m-proof}

We probe this matter numerically, by plotting a few examples with the Runge-Kutta method (using MAPLE).

\begin{longtable}{ccc} 
\caption{Examples of plots using the Runge-Kutta method}
\\ \hline 
\small{values of parameters}&\small{graphs of envelopes}&\small{graph of error } \\ 
&\small$y_{RK+}\,$(dots), $y_{RK-}\,$(line) &\small$err_y(x)=y_{RK+}-y_{RK-}$
 \\ \hline &&\\ 
$\begin{array}{l}\bb=10,\\ \tu=0.3,\\ \BB=0.73,\\ max.err_y=0.035,\\ \frac{\ln(max.err_y)}{\BB}=-4.55.
\end{array}$
&
\begin{minipage}[c]{0.275\linewidth}
\centering
{\includegraphics[width=0.75\textwidth]{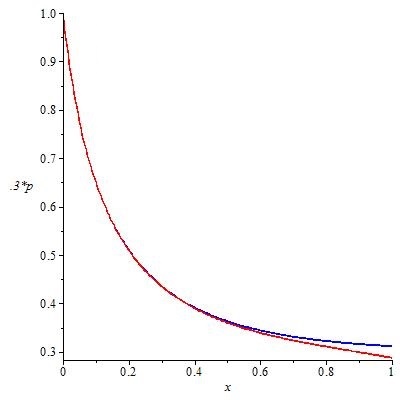}}
\end{minipage}
&
\begin{minipage}[c]{0.275\linewidth}
\centering
{\includegraphics[width=0.75\textwidth]{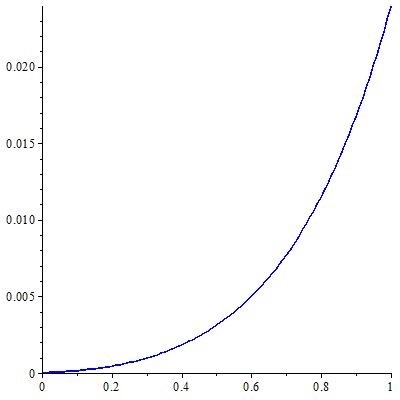}}
\end{minipage}
 \\ \hline && \\
$\begin{array}{l}\bb=55,\\ \tu=0.1,\\ \BB=0.78,\\ max.err_y=0.038,\\ \frac{\ln(max.err_y)}{\BB}=-4.18.
\end{array}$
&
\begin{minipage}[c]{0.275\linewidth}
\centering
{\includegraphics[width=0.75\textwidth]{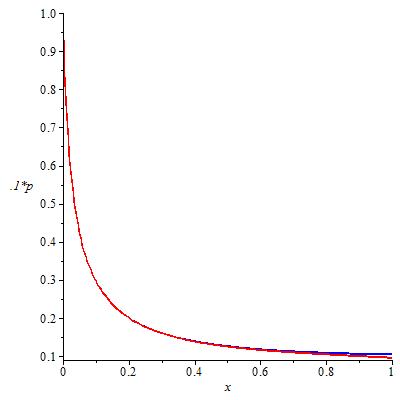}}
\end{minipage}
&
\begin{minipage}[c]{0.275\linewidth}
\centering
{\includegraphics[width=0.75\textwidth]{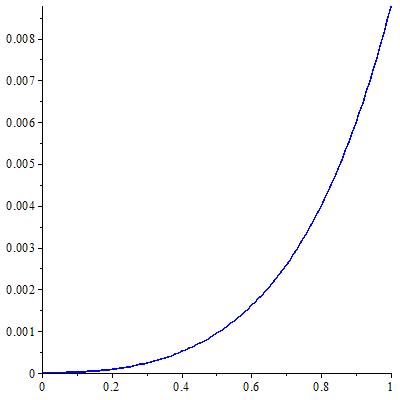}}
\end{minipage}
 \\ \hline && \\
$\begin{array}{l}\bb=30,\\ \tu=0.7,\\ \BB=7.85,\\ max.err_y=4.7\cdot 10^{-11},\\ \frac{\ln(max.err_y)}{\BB}=-3.02.
\end{array}$
&
\begin{minipage}[c]{0.275\linewidth}
\centering
{\includegraphics[width=0.75\textwidth]{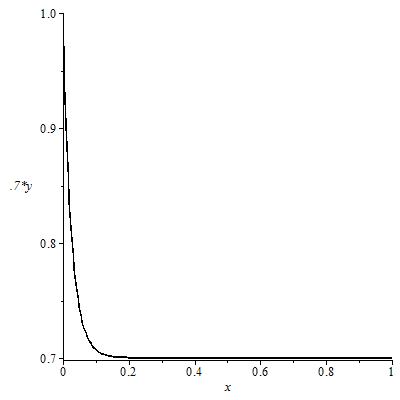}}
\end{minipage}
&
\begin{minipage}[c]{0.275\linewidth}
\centering
{\includegraphics[width=0.75\textwidth]{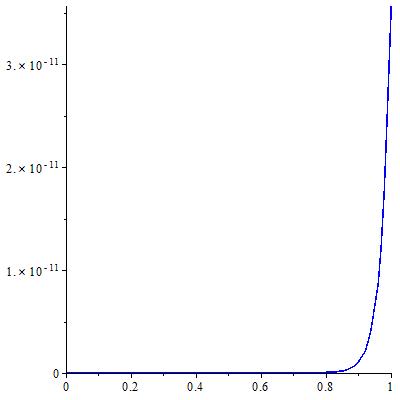}}
\end{minipage} \\ \hline 
\label{tab:err}
\end{longtable}

The conclusion is that the Runge-Kutta method can be successfully run, as long as there is enough computational power available; high values of the parameters require long computational time. It's certainly remarkable that the error decreases exponentially fast with $\BB$. (The last graphs of the envelopes overlap.)


\subsection{The boundary layer condition}

The estimates~\eqref{eq:dery1} are essential for understanding the behaviour of the solution to the equation~\eqref{eq:y-het}, especially the boundary layer phenomenon.  

\begin{tabular}{ll}
\hspace{-3ex}
\begin{minipage}{.6\textwidth}
Note that the secant line between $(0,1)$ and $(1,\TT)$ is above the graph of $y_\het$, so $y_\het'(1)>\TT-1$. Now consider the angles, denoted $a$ and $s$ in the figure. For a parameter $\kappa>1$, the {\em $\kappa$-boundary layer} property should be that 
\begin{m-eqn}{
a<\frac{s}{\sqrt{2}\kappa}=\frac{\tan^{-1}\bigl(\frac{1}{T-1}\bigr)}{\sqrt{2}\kappa}.
}\label{eq:geom-k-bl}
\end{m-eqn}
\end{minipage}
&
\begin{minipage}{.295\textwidth}
\includegraphics[width=.85\textwidth]{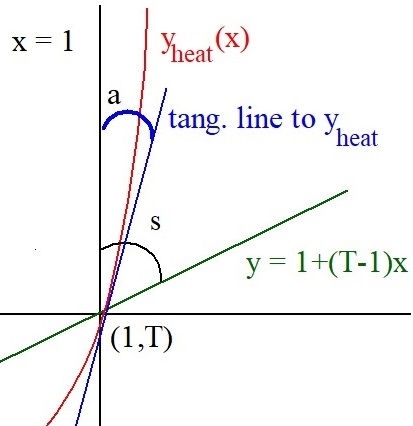}
\end{minipage}
\end{tabular}

For $\TT$ ---same for $\tu=\TT^{-1}$--- close to $1$, the function $y_\het$ is approximately constant $1$. Thus, we are interested in larger values, such as $T\ges 1.5$ that is $\tu\les 2/3$. Since $y_\het'(1)=a^{-1}$ is approximately $\sqrt{2}BT^{2.5}$, we replace the inequality above with $y_\het'(1)>\sqrt{2}\kappa\TT$, which yields: 
\begin{m-eqn}{
BT^{1.5}>\kappa
\quad\Leftrightarrow\quad 
b>\sqrt{5}\kappa.
}\label{eq:k-bl}
\end{m-eqn}
We choose the value $\kappa=35$ for our parameter: then~\eqref{eq:geom-k-bl} implies, for $T>2$, that the angle $a$ is less than $1^\circ$. From a visual standpoint, at least for $x$ very close to $1$, this should be steep enough to claim the existence of a ``boundary layer''. 

Now, let us compare this geometric approach with the analytic data: we plot the graph of $\typ$ for several values of the parameters $\BB$ and $\TT$. Since $\typ$ is greater than $y_\het$, the graph of latter is below the former, thus $y_\het$ is ``at least as horizontal'' about $x=0$ and ``at least as vertical" about $x=1$. We verify whether, from a visual viewpoint, the boundary layer property is imposed only by the size of the derivative at $x=1$. The will see that we must slightly adjust our previous considerations. \smallskip

\begin{tabular}{ll}
\hspace{-3ex}
\begin{minipage}{.625\textwidth}
The estimate for $y_\het'(0)$ shows that already for low values of $\BB$ (e.g. $\BB\ges0.5$), the graph of $y_\het$ starts off almost horizontally at $x=0$. This is not true for low values of $\BB$\/! 
\end{minipage}
&
\begin{minipage}{.27\textwidth}
\includegraphics[width=.85\textwidth]{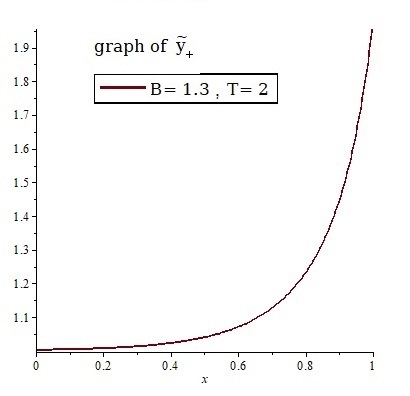}
\end{minipage}
\end{tabular}

\begin{tabular}{ll}
\hspace{-3ex}
\begin{minipage}{.625\textwidth}
On the right, $\BB\TT^{1.5}>35$ is large, yet there is no (visual) boundary layer. (The plot is only over the interval $x\in[0.9,1]$). The explanation for this matter is the following: {\em the condition $y_\het'(1)\gg0$ doesn't necessarily imply that the graph ``falls vertically'' along the line $x=1$}.\\ 
The idea which emerges is that, to ensure steep-falling, one should impose that $y_\het'$ is large at a suitable value $x<1$. (Note that, since $y_\het'$ is increasing, $y_\het'(1)$ is large, too.)
\end{minipage}
&
\begin{minipage}{.27\textwidth}
\includegraphics[width=.85\textwidth]{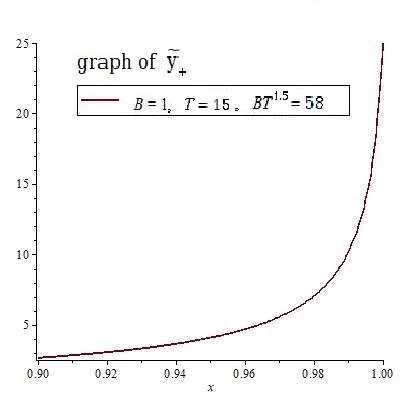}
\end{minipage}
\end{tabular}\medskip 

\begin{m-definition}\label{def:bl}
We say that $y_\het$ possesses the {\em boundary layer property} if the following condition holds:
$$
\BB\cdot T^{1.25}\ges 50\quad\text{or}\quad b\sqrt[4]{t}>110. 
$$
\end{m-definition}

\begin{itemize}
\item Let $h:=y_\het^{-1}(\sqrt{\TT})<1$. The estimate~\eqref{eq:dery1} shows that our definition  basically means 
$$y_\het'(h)\ges50\sqrt{2}>\tan(89^\circ)\approx 57.29,$$ 
so the angle $a$ in~\eqref{eq:geom-k-bl} is less than $1^\circ$, too. 

\item Our benchmark for declaring the ``boundary layer" was the angle $1^\circ$.  One obtains steeper boundary layer ---the graph of $y_\het$ is even more vertical--- by increasing the lower bound of $\BB\TT^{1.25}$. 
\end{itemize}\medskip

We conclude with another geometric description of the boundary-layer property. 

\begin{m-corollary}
The existence of the boundary layer implies that the function dramatically decreases over a very short interval. The ratio between the variation of $y_\het$ over the interval $[1-\sqrt{5}/b,1]$ and its total variation over $[0,1]$ is at least $0.5$. 
\end{m-corollary}

\begin{m-proof}
We observe that $\xi=\sqrt{5}/b=t^{1.5}/B\les t^{0.25}/50\les0.02$, so the indicated interval is indeed short. Moreover, since $\typ\ges y_\het$ and $\tBp\ges\BB$, we have: 
$$
\frac{\TT-y_\het(1-\xi)}{\TT-1}
\ges
\frac{\TT-\typ(1-\xi)}{\TT-1}
\ges 
\frac{1-t\,\cdot\,{\tanh(\, \tanh^{-1}(\tu^{1.5})+1.5\sqrt{2}t^{1.5})\,)^{-\frac{2}{3}}}}{1-\tu}>0.5.
$$
The last inequality is obtained by plotting the function (for $\tu\in[0,1]$) and reading off its minimum ($\approx0.53$).
\end{m-proof}\medskip


\subsection{Inconveniences}\label{ssct:inconv}

We justify the necessity for further investigation. 
\begin{itemize}
\item 
First, one can't hope to obtain the opposite inequality~\eqref{eq:tBB} ---it's clear for $T\gg0$---, so the IVP~\eqref{eq:odey} is not suitable for deducing a lower bound to $y_\het$.
\item Second, the sensitive dependence of the problem on the parameters $b,\tu$, makes that $\typ$ above is typically a rather loose estimate for $y_\het$. \\ 
In Table~\ref{tab:y+res}, we plotted $\typ$ and its residue 
$$
\resd(\typ):={\tBp}^2(5\typ^4-4\typ-\typ^{-2})
-B^2(5\typ^4-5),
$$
for $\BB=13,\TT=3$.  
We are in the boundary layer case $\BB\TT^{1.25}\approx51$, the graph falls vertically, but the residue is large near $x=1$; its $L^2$-norm is $\|\resd(\typ)\|_{L^2{([0,1]})}\approx 36.7$.
\end{itemize}
\smallskip

\begin{longtable}{ll}
\caption{Graphs of $\typ$ and $\resd(\typ)$.}
\\[-1ex] 
\begin{minipage}{.275\textwidth}
\includegraphics[width=.75\textwidth]{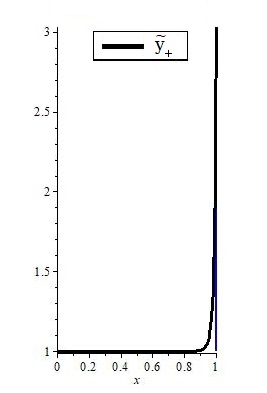}
\end{minipage}
&
\begin{minipage}{.275\textwidth}
\includegraphics[width=.9\textwidth]{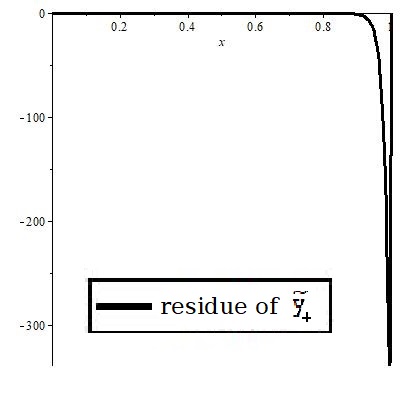}
\end{minipage}
\label{tab:y+res}
\end{longtable}


\section{The heat conduction equation: II}\label{sct:heat2}

Now we refine the previous methods.


\subsection{A boundary value problem\,II}\label{ssct:bvp2} 

The change of variables $z=y^{1.5}$ in~\eqref{eq:odey} yields the ODE $2z'=3\sqrt{2}\BB(z^2-\CC^2)$, which requires integrating 
$$
\frac{2\rd z}{z^2-C^2}=\frac{\rd z}{C}\Big[ \frac{1}{z-C} - \frac{1}{z+C} \Big].
$$ 
Here we consider a one-parameter deformation of this $1$-form, namely:
$$
\frac{\rd z}{C}\Big[ \frac{1}{z-C\ee} - \frac{1}{z+C} \Big]=\frac{(1+\ee)\rd z}{z^2+(1-\ee)Cz-\ee C^2},
$$ 
where $\ee<1$ is introduced to shift the pole of the initial $1$-form away from $z=C\approx 1$, since poles yield major computational errors. The function $z$ solves 
$$
\frac{\rd z}{C}\Big[ \frac{1}{z-C\ee} - \frac{1}{z+C} \Big]=3\sqrt{2}\BB\rd x,
\quad z(0)=z_0\approx 1,\;z(1)=\TT^{1.5},
$$
so it satisfies the equation 
\begin{m-eqn}{
\frac{z-\ee\CC}{z+\CC}
=\frac{\TT^{1.5}-\ee\CC}{\TT^{1.5}+\CC}\cdot e^{-3\sqrt{2}\BB\CC(1-x)},
}\label{eq:z-ee}
\end{m-eqn}
with $\CC$ determined by the equality:
\begin{m-eqn}{
\frac{\TT^{1.5}-\ee\CC}{\TT^{1.5}+\CC}
=\frac{z_0-\ee\CC}{z_0+\CC}\cdot e^{2\LL\CC},\quad\LL:=1.5\sqrt{2}\BB.
}\label{eq:C-ee}
\end{m-eqn}
We recover the situation described in Remark~\ref{rmk:cc}.
 
For the initial function $y=z^{2/3}$, we are considering the solution to the new BVP,   
\begin{m-eqn}{
\scalebox{.95}{$
\frac{y'}{\sqrt{2}}=\frac{2\tBB}{1+\ee}\bigl(y^{2.5}+(1-\ee)\CC y-\ee\CC^2y^{-0.5}\bigr),\quad y(0)=1,\;y(1)=\TT, 
$}
}\label{eq:odee1}
\end{m-eqn}
whose second derivative satisfies:
\begin{m-eqn}{
\scalebox{.95}{$
y''=\frac{4\tBB^2}{(1+\ee)^2}\Bigl( 
5y^4+7(1-\ee)\CC y^{2.5}-2(4\ee-1-\ee^2)\CC^2y-\ee(1-\ee)\CC^3y^{-0.5}-\ee^2\CC^4y^{-2}
\Bigr).
$}
}\label{eq:odee2}
\end{m-eqn}
The discussion in Section~\ref{ssct:bvp1} corresponds to $\ee=C=1$. We introduce the notation: 
$$
\BB=q\cdot\tBB. 
$$
The residue function of $y$ satisfies the formula:
$$
\begin{array}{rl}
\frac{\resd(y)}{\BB^2}&=R(y,q,\CC,\ee)
=
\bigl(1-q^2\frac{(1+\ee)^2}{4}\bigr)\cdot5y^4+7(1-\ee)\CC y^{2.5} -2(4\ee-1-\ee^2)\CC^2y
\\[2ex]&
-\ee(1-\ee)\CC^3y^{-0.5}-\ee^2\CC^4y^{-2}+5q^2\frac{(1+\ee)^2}{4},
\end{array}
$$
involving the parameters $q,\ee,\CC$, which must be finely tuned so that Proposition~\ref{prop:+-} will be eventually satisfied: the solution to the BVP~\eqref{eq:odee2} is a lower (resp. upper) bound for the solution of~\eqref{eq:y-het} if $R(y,q,\CC,\ee)$ is positive (resp. negative) for $y\in[0,1]$. 

\begin{m-procedure}\label{proc:qec}
\begin{enumerate}
\item The value of $q$ is determined by the derivative condition  $$\frac{y'(1)}{\sqrt{2}}\approx\BB\sqrt{T^5-4T+(4+\dta)},$$  
and by ensuring the correct inequality; the bounds~\eqref{eq:dery1} for $\dta$ are essential.  
\begin{itemize}
\item for $\les$, we impose $y'(1)/\sqrt{2}\les q\tBB\cdot\sqrt{T^5-4T+4}$.
\item for $\ges$, we impose $y'(y)/\sqrt{2}\ges q\tBB\cdot\sqrt{y^5-4y+(4+\Delta)}$, for $y=T,\sqrt{T}$;
\end{itemize}

\item The deformation parameter $\ee$ will be typically close to $1$, and it's determined by the curvature/residue condition $R(y,\dots)=r\BB^2,$ where $r$ is a small real number (positive, resp. negative). 

\item Finally, $\CC$ is determined by the boundary condition $y(0)\approx 1$ (cf. Section~\ref{sct:back});  one expects it to be close to $1$, too. Note that $\CC$ is an explicit function of $q,\ee$.
\end{enumerate}
\end{m-procedure}


\subsection{Analysis of the residue}\label{sssct:sign-res} 

The sign of the residue function is a decisive matter. In explicit situations, the simplest way to settle the issue is by plotting the function. Here we analyse the sign of this expression (dependence on parameters) by relying on analytical tools. The polynomial in four variables above has the following properties:

\begin{itemize}
\item $R(1,1,1,1)=0,\;\frac{\rd R}{\rd C}(1,1,1,1)<0,\; \frac{\rd R}{\rd \ee}(1,1,1,1)<0;$
\item $R(1,1,\CC,\ee)= \CC^2(2+\CC-\CC^2)\ee^2-\CC(7+8\CC+\CC^2)\ee+5+7\CC+2\CC^2;$ 
\item The quadratic equation $R(1,1,1,\ee)=0$ admits two roots: $1$ and $7$. For small variations of the parameters, these roots will persist, and we denote by $\ee(y,q,\CC,r)$ the root of $R(y,q,\CC,\ee)=r$ which is close to $1$. For shorthand, we abusively denote it by $\ee(y)$, since we are mainly interested in the dependence in $y$.
\end{itemize}

\begin{m-lemma}\label{lm:es}
We denote $s:=y^{-1}\in[0,1]$, and let: 
$$
\tilde\ee_-(s):=1-\frac{4}{5}s^{3}+\frac{3}{5}s^{4},\qquad
\tilde\ee_+(s):=1-\frac{4}{5}s^{3}+\frac{8-3s}{5}s^{4}.
$$
Then the following statements hold: 
\begin{itemize}
\item $R(s^{-1},1,1,\tilde\ee_-(s)), R(s^{-1},1,1,\tilde\ee_+(s))$ are continuous, for $s\in[0,1]$.
\item They satisfy the inequalities:
$$
\forall\,s\in[0,1],\qquad R(s^{-1},1,1,\tilde\ee_-(s))\ges2,\quad R(s^{-1},1,1,\tilde\ee_+(s))\les-\frac{2}{3}.
$$
\end{itemize}
Thus, for small values of the parameters $q,\CC,r$, the root $\ee(s)$ of $R(s^{-1},q,\CC,\ee)=r$ belongs to the interval $[\tilde\ee_-(s),\tilde\ee_+(s)]$. 
\end{m-lemma}
One should still clarify what's the meaning of `small' in the statement. The series expansion of the root of $R(s^{-1},1,\CC,r)=0$ is $\ee(y)=1-\frac{4}{5}\CC^2s^3+\Big(1-\frac{r}{5}\Big)s^4+O(s^5)$. One may check that the Lemma still holds for $1-\CC\les t^3/10$ and $|r|\les0.5$. 

\begin{m-proof}
The continuity reduces to the fact that the series development of $R$ contains only positive powers of $s$. The indicated inequalities are obtained by plotting the corresponding graphs. 

For the second statement, a computation shows that the power series expansion of the root of $R(s^{-1},1,\CC,r)=0$ is $\ee(y)=1-\frac{4}{5}\CC^2s^3+\Big(1-\frac{r}{5}\Big)s^4+O(s^5)$. 
\end{m-proof}

Lemma~\ref{lm:es} is important to understand the monotonicity of the residue function.  Below we compute the components of the gradient of $R$; the abbreviation $\lot$ stands for `lower order terms'. Also, we think off $C$ as function of $q$ and $\ee$, since $C$ is determined by the boundary condition at $x=0$ (see Procedure~\ref{proc:qec}); a {\em very rough} approximation would be $\tanh(1.5\sqrt{2}Bq^{-1})$.
$$
\scalebox{.85}{$
\begin{array}{ll}
\disp\frac{\rd R}{\rd y}&
= 
\Big( 1-q^2\frac{(1+\ee)^2}{4} \Big)\cdot 20y^3
+17.5(1-\ee)\CC y^{1.5} - 2(4\ee-1-\ee^2)\CC^2+\lot;
\\[2ex] 
\disp{\Bigl.\frac{\rd R}{\rd q}\Bigr|}_{y=T}&
=-10q\frac{(1+ee)^2}{4}T^4+7(1-\ee)\frac{\rd C}{\rd q}T^{2.5}+\lot\approx -10T^4-14(1-\ee)BT^{2.5};
\\[2ex] 
\disp{\Bigl.\frac{\rd R}{\rd \ee}\Bigr|}_{y=T}&
=-5q^2\frac{1+\ee}{2} T^4-7CT^{2.5}+7(1-\ee)\frac{\rd C}{\rd \ee}T^{2.5}+\lot\approx -5T^4.
\end{array}
$}
$$
For $1-\ee\approx\frac{4}{5}t^3$, as suggested by the Lemma, $\frac{\rd R}{\rd y}$ is mostly positive, so $y\mt \resd(y)$ will be increasing almost over the whole $[0,\TT]$. (Here we think off $y$ as a variable, rather than a function of $x$.) Thus, for $y=y(x)$ an approximate solution of the heat equation, $\resd(y(x))$ increases dramatically, for $x$ near $1$.  

Furthermore, the dominant component of the gradient is $\frac{\rd R}{\rd q}$; the factor of $T^4$ is responsible for large variations of the residue for tiny variations of the parameters. The only way to ascertain the increase/decrease of $R$ is by changing $q,\ee$ along the gradient flow.


\section{Bounding the exact solution $y_\het$}\label{sct:bounds}

\subsection{A minorant $y_-$ of $y_\het$}\label{ssct:ee-lower}

\nit\underbar{(i)~determine $q$}\; 
In the sequel, we let  
$$
M(y^{-1}):=\sqrt{1-5y^{-4}+(4+\Delta)y^{-5}},\; q(y,\ee,\CC):=\frac{2(1+(1-\veps)\CC\cdot y^{-1.5}-\veps\CC^2\cdot y^{-3})}{(1+\veps)\cdot M(y^{-1})}.
$$
The slope inequality condition $y_\het(y)\les\frac{y'(y)}{\sqrt{2}}$ is satisfied as soon as $q\les q(y,\ee,\CC)$; the right-hand side is decreasing in $\CC$, so we may replace it by its upper bound $\veps^{-1}$ (see equation~\eqref{eq:CC-bounds}). Hence the slope inequality at $y$ is satisfied for 
$$
q\les q(y,\ee,\ee^{-1}).
$$

\subsubsection{A global minorant}\label{sssct:global-min}

For $\veps=1,q=1$, the residue is negative; for $\veps=0,q=1$, the residue is positive. We need values of $\veps,q$ in between, for which the residue positive, for all  $y\in[0,1]$, see Proposition~\ref{prop:+-}. 

\begin{m-proposition}\label{prop:global-min}
Let $b>0, t\in(0,1)$,  $\tilde\ee_{-}:=0.75$. Define
$$
\tilde q_{-}:=\min\Bigl(1.1, q(T,\tilde\ee_-,\tilde\ee_-^{-1})\Bigr),
\;
\tBB_-:=\tilde q_-^{-1}\BB .
$$
We define $\tilde\CC_-:=\CC_-(\tilde q_-,\ee_-)$ as in Section~\ref{sct:back}. Then the function $\tilde y_-:=y_{\tBB_-,\CC_-,\veps_-}$ defined by equation~\eqref{eq:ODE} is a global lower envelope for $y_\het$.
\end{m-proposition}

\begin{m-proof}
We prove that the residue function is positive that is, $R(y,q_-,C,\veps_-)>0.$ Since $R$ is decreasing in $C$, we may replace it by an upper bound; by~\eqref{eq:CC}, we have $C\les\veps^{-1}$. By the previous discussion, $q$ must satisfy the inequality: 
$$
q\les q(T,\veps_-,\veps_-^{-1})\les
\frac{2(1+(1-\veps_{-})C\cdot t^{1.5}-\veps_{-}C^2\cdot t^3)}{(1+\veps_{-})\cdot M(t)}.
$$
Since $\tilde q_-\les1.1$, too, the residue condition becomes $R(y,1.1,\tilde\ee_-^{-1},\tilde\ee_-)\ges0$, for $y\in[1,T]$. The only negative trouble-making term is  $2(4\tilde\ee_--1-\tilde\ee_-^2)\tilde\ee_-^{-1}y=\frac{23}{6}y$, and it's dominated either by the first two terms (for $y>1.2$) or by the last one (for $y\in[1,1.2]$). Alternatively, one may test the inequality graphically.
\end{m-proof}


\subsubsection{A partial minorant}\label{sssct:part-min}

When attempting to improve the precision of the lower bound, one inevitably faces the strong increasing tendency of $y\mt R(y,\dots)$, as $y$ approaches $T$, which practically makes impossible to have a global lower bound with small residue. 

The way out is to determine a partial lower bound along the boundary layer $\sqrt{T}\les y\les T$, where is concentrated the interesting information. We implement this strategy by imposing positive curvature and slope conditions at $y=\sqrt{T}$. They yield the values of $\ee$ and $q$, respectively; since $R$ is increasing in $y$, the same conditions will hold along $[\sqrt{T},T]$. Concretely, we proceed as follows: 
\begin{itemize}
\item Our relevant functions $q,R$ are decreasing in $\ee,\CC$, so we fix the `worst' maximal possible values: $\ee_0=\tilde\ee_{+}(t),\;\CC_0:={\tilde\ee_-(t)}^{-1}$. They are needed for writing the slope and curvature conditions at $y=\sqrt{T}$:
$$
\ee_1:=\ee_0,\;q_1:=q(\sqrt{T},\ee_0,\CC_0),\;\CC_1:=\CC_-(\ee_1,q_1).
$$
Recall from Lemma~\ref{lm:es} that $R$ is negative for $\ee=\ee_1, y=T$ ---thus also for $y=\sqrt{T}$---, so the true value, making $R$ positive, is very slightly less than that. Therefore the slope condition, imposed by the value of $q_1$, will be satisfied by the true value, too. 
\item The precision increases if $(\ee_1,q_1)$ is modified in the direction of the gradient of $R$.
\end{itemize}

\begin{m-proposition}\label{prop:partial-}
Let $(\ee_-,q_-,\CC_-):=(\ee_1,q_1,\CC_1)$, and $\BB_-:=q_-^{-1}\BB$. The function 
$$
y_-(x):=\biggl[ 
\CC_-\cdot\frac{\TT^{1.5}\cdot (1+\ee_- e^{3\sqrt{2}\BB_-\CC_-(1-x)})+\ee_- \CC_-(e^{3\sqrt{2}\BB_-\CC_-(1-x)}-1)}{T^{1.5}(e^{3\sqrt{2}\BB_-\CC_-(1-x)}-1)
+\CC_-(e^{3\sqrt{2}\BB_-\CC_-(1-x)}+\ee_-)}
 \biggr]^{2/3}
$$
satisfies the differential equation~\eqref{eq:odee1}, and is a partial lower envelope of the exact solution $y_\het$, for values $y\in[\sqrt{T},T]$.  
\end{m-proposition}

\begin{m-proof}
We apply Proposition~\ref{prop:+-} and the discussion following it.
\end{m-proof}


\subsection{An majorant $y_+$ for $y_\het$}\label{ssct:ee-upper} 

A global upper envelope was already constructed in Proposition~\ref{prop:dery1}.  

\subsubsection{A partial majorant}\label{sssct:part-max}

We construct a partial  upper envelope along the boundary layer. The monotonicity of $y\to R(y,\dots)$ implies that the residue function is automatically negative about $x=1$, as soon as $\resd(y=T)\les0$.  

We repeat the previous steps, but with the opposite inequalities; this task will be more delicate than before. Recall that, for $q=\CC=\ee=1$, the function $\typ$ at Section~\ref{ssct:bvp1} is an upper bound for $y_\het$; its major shortcoming was the (negative) magnitude of the residue about $x=1$. To fix this issue, we need more accurate derivative and curvature conditions at $x=1$, by following the gradient of $R$.  To implement this desideratum, we need a starting values for $q,\ee,C$. The function $\typ$ is defined  (cf. equation~\eqref{eq:tBp}) using the factor 
$$
\tilde q_+=\Bigl[\frac{(1-\tu^{3})(5+\tu^{3})}{5(1-\tu^4)}\Bigr]^{1/2}=1-\frac{2}{5}t^3+\frac{1}{2}t^4-\dots
$$
For this reason, we take 
$$
q_{0+}:=1-\frac{3}{5}t^3
$$
as the `worst', smallest possible $q$-factor, which must be increased. For the `worst' $\ee$-value, Lemma~\ref{lm:es} readily implies that 
$$
\ee_{0+}:=1-\frac{4}{5}t^3
$$
is quite an accurate choice. Finally, we need an (absolute) lower bound for $C$, and this is $\CC_{0+}:=c_0(L)$. Note that $c_0(L)$ is defined as an infinite iteration of the $\tanh$-function; for numerical calculations, we'll define $\CC_{0+}$ as $30$ iterations of $\tanh(0.99L)$:
$$
C_{0+}:=\tanh(0.99L\cdot\tanh(0.99L\cdot(\dots\cdot\tanh(0.99L)\dots))).
$$

\medskip 
\nit\underbar{(i)~determine $\ee$}\;  We let 
$$
\ee_+:=\ee(\TT,q_{0+},C_{0+},0)
$$ 
that is, $\ee_+$ satisfies $R(\TT,q_{0+},C_{0+},\ee_+)=0$. Since $q_{0+}$ is the smallest possible $q$-factor and $R$ is decreasing in $q$, we'll have $R(\TT,q,C_{0+},\ee_+)<0$ for our subsequently computed values of $q$. 

\medskip
\nit\underbar{(ii)~determine $q$}\; The derivative condition at $x=1$ is $\frac{y'(1)}{\sqrt{2}}\les\BB\sqrt{\TT^5-5\TT+4}<\frac{y_\het'(1)}{\sqrt{2}},$ so we must have  
$$
q\ges q_+(\ee,\CC):=
\frac{2}{1+\ee}\cdot\frac{\TT^{2.5}+(1-\ee)\CC\TT-\ee\CC^2\TT^{-0.5}}{\sqrt{\TT^5-5\TT+4}}.
$$
Since the right-hand side is decreasing with $\ee,\CC$, we define 
$$
q_{1+}:=q_+(\ee_{0+},\CC_{0+}). 
$$
This choice  ensures that the approximate function is greater than $y_\het$, about $x=1$. We remark that it's not possible to neglect $\CC,\ee$; doing so leads to severely wrong calculations. 

\medskip 
\nit\underbar{(iii)~determine $\CC$}\; 
The last condition we impose is $y_+(0)\ges 1$. The formula~\eqref{eq:C-ee} shows that this boils down to: 
$\frac{\TT^{1.5}-\ee\CC}{\TT^{1.5}+\CC}\ges \frac{1-\ee\CC}{1+\CC}\cdot e^{2\LL\CC}.$ Our discussion in Section~\ref{ssct:treq-2} yields: 
\begin{m-eqn}{
\begin{array}{l}
\LL_{1+}:=q_{1+}^{-1}\LL,\quad\zeta_{1+}:=\frac{\TT^{1.5}+\frac{1-\ee_{1+}}{2\ee_{1+}}}{1+\frac{1-\ee_{1+}}{2\ee_{1+}}},\quad 
\\[2ex] 
\cc_{1+}:=\cc_0(\zeta_{1+},\ee_{1+}^{-1}\LL_{1+}),\quad 
\CC_{1+}:=\CC_+(\ee_{1+},\cc_{1+}).
\end{array}
}\label{eq:C+}
\end{m-eqn}

\medskip
\nit\underbar{(iv)~iterations}\; The estimates greatly improve by iterating the process. The value $\ee_{+}$ above is quite precise, so we keep it. One improves $q$ by using the gradient of $R$,
$$
q_{2+}:=q_{1+}+\frac{R(\TT,q_{1+},\CC_{1+},\ee_+)}{20\cdot B^2}t^4,
$$
and then $C_{2+}$ is determined by $\ee_+,q_{2+}$, as in~\eqref{eq:C+} above. The strange-looking modification for $q$ is motivated by the estimates for $\frac{\rd R}{\rd q}$.

\begin{m-proposition}\label{prop:partial+}
Let $\ee_+$ be as above and denote by $q_+$ and $\CC_+$ the values obtained by two iterations; let $\BB_+:=q_+^{-1}B$. The function 
$$
y_+(x):=\biggl[ 
\CC_+\cdot\frac{\TT^{1.5}\cdot (1+\ee_+ e^{3\sqrt{2}\BB_+\CC_+(1-x)})+\ee_+ \CC_+(e^{3\sqrt{2}\BB_+\CC_+(1-x)}-1)}{T^{1.5}(e^{3\sqrt{2}\BB_+\CC_+(1-x)}-1)
+\CC_+(e^{3\sqrt{2}\BB_+\CC_+(1-x)}+\ee_+)}
 \biggr]^{2/3}
$$
is a partial upper bound of the exact solution $y_{heat}$, for values $y\in[\sqrt{T},T]$.
\end{m-proposition}

\begin{proof}
We apply Proposition~\ref{prop:part+-} and the discussion following it.
\end{proof}


\section{Back to the heat conduction equation}\label{sct:back-heat}

This study is motivated by the  heat conduction equation~\eqref{eq:uhet}. The bounds for $y_\het$ immediately  yield bounds for $u_\het$: we just set $u=t\cdot y(1-x)$ everywhere. Thus global/partial upper/lower envelopes for $y_\het$ determine, respectively, the same type of envelopes for $u_\het$. The error for estimating $u_\het$ is $t$ times the error for $y_\het$, thus smaller. The maximal computational errors for $y_\het, u_\het$ are compared in Table~\ref{tab:maxerr-yu}, for various values of $b,t$. The Runge-Kutta method  (see Section~\ref{ssct:shoot}) has precision $O(t\cdot e^{-3B})$. 

\subsection{Case {\it t }{\rm = 0}} 

An issue which naturally raises is to determine an approximate solution of the initial heat equation, for the limiting value $t=0$. In physical sense, this situation occurs when the ambient (final) temperature is negligible compared to the temperature of the source. Note that $B=O(t^{1.5})$, it becomes small for $t$ approaching $0$, so this situation is exactly opposite to the `$B$ large' case.

\begin{m-proposition}
An upper envelope for the BVP, 
$$\;u_0''=b^2u_0^4,\; u_0(0)=1,\;u_0(1)=0,\;$$ is 
$$\disp \tilde u_{0+}(x):={\Big[1+\frac{1.5\sqrt{2}}{\sqrt{5}}bx\Big]}^{-2/3}.
$$
The largest deviation from the exact solution of the BVP ---the maximal error--- is attained at $x=0$, and  decreases with $b$: 
$$
max.err_{\tilde u_{0+}}={\Big[1+\frac{1.5\sqrt{2}}{\sqrt{5}}b\Big]}^{-2/3}
\approx{\Big[1+0.95\cdot b\Big]}^{-2/3}.
$$
\end{m-proposition}
Numerical data in the next section shows that the global upper bound $\tup$ for $y_\het$ is quite precise already. Below, we'll see a convenient method to determine accurate lower envelope, too.

\begin{m-proof}
Just compute $\disp\lim_{t\to 0^+}t\cdot\typ(1-x)$. For the last statement, $u_0$ satisfies the IVP  
$$
u_0'=-\frac{\sqrt{2}b}{\sqrt{5}}\cdot\sqrt{u_0^5+\gamma}, 
\quad u_0(0)=1, u_0(1)=0,
$$ 
with $\gamma>0$ is determined by the boundary conditions, and 
$\tilde u_{0+}$ satisfies $$u_{0+}'=-\frac{\sqrt{2}b}{\sqrt{5}}u_{0+}^{5/2},\;u_{0+}(0)=1.$$ Note that $\frac{\sqrt{2}b}{\sqrt{5}}\cdot\sqrt{u^5+\gamma}-\frac{\sqrt{2}b}{\sqrt{5}}u^{5/2}$ is decreasing with $u$, hence the difference of slopes is the smallest at the $x=0$ end and increases towards $x=1$ end of $[0,1]$. 
\end{m-proof}

There is yet another way to numerically determine a better approximation for $u_0$, by using the Runge-Kutta method and the estimates~\eqref{eq:dery1}. For $t>0$ we have $u_\het(x)=t\cdot y_\het(1-x)$, so we approximate  ${\gamma}$ above by 
\begin{m-eqn}{
{\Gamma}:=\lim_{t\to0^+}\biggl(t\cdot\frac{\sqrt{5}}{\sqrt{2}b}\typ'(0)\biggr)^2
={\biggl(  1+\frac{1.5\sqrt{2}\,b}{\sqrt{5}}  \biggr)}^{-10/3}
=
\biggl[\frac{\sqrt{5}}{\sqrt{2}\,b}\cdot(-\tilde u_{0+}'(1))\biggr]^2.
}\label{eq:Gamma}
\end{m-eqn}
The numerical function $w_0$ defined by the IVP  
\begin{m-eqn}{
w_{0,\Gamma}'=-\frac{\sqrt{2}\,b}{\sqrt{5}}\sqrt{w_{0,\Gamma}^5+\Gamma},\quad w_{0,\Gamma}(0)=1, 
}\label{eq:w0}
\end{m-eqn}
satisfies the ODE $w_{0,\Gamma}''=b^2 w_{0,\Gamma}^4$, the same as $u_0$. It approximates the exact solution $u_0$ at least as well as $\tilde u_{0+}$ does, and can be plotted by using the Runge-Kutta method.  
 
\begin{longtable}{cc}
\caption{Numerical approximations of $u_0$ and their errors.}\\ 
\kern-1.5ex\scalebox{0.95}{
\begin{tabular}{|r|c|c|}\hline 
$b$&$w_{0,\Gamma}(1)$&$\tilde u_{0+}(1)-w_{0,\Gamma}(1)$
\\  \hline 
10 &   0.193  & 0.015  
\\ \hline 
30 & 0.097   & 0.008
\\ \hline 
70 & 0.055  &  0.004
\\ \hline 
100 & 0.044 & 0.003
\\ \hline 
500 & 0.015 & 0.001
\\ \hline 
1000 & 9.5 E-3 & 7.8 E-4
\\ \hline 
50000 & 7.0 E-4& 5.7 E-5 
\\ \hline 
100000 & 4.4 E-4 & 3.6 E-5  
\\ \hline 
\end{tabular}
}
&\kern-2.5ex
\begin{minipage}[c]{0.55\linewidth}
{\includegraphics[width=0.95\textwidth]{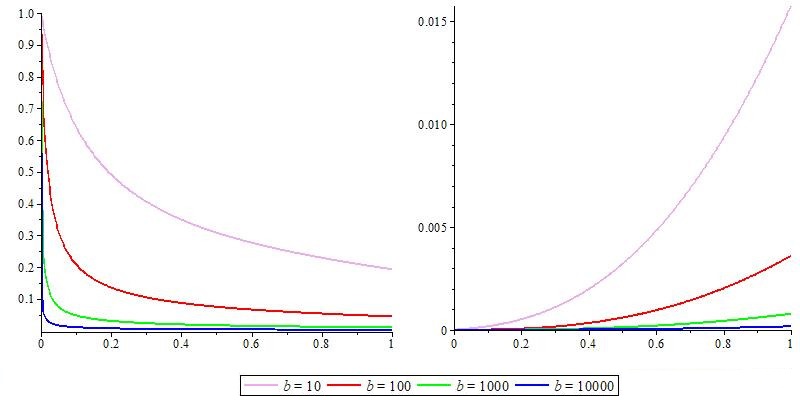}}
\end{minipage}
\\   \renewcommand{\arraystretch}{1.5}
{\scriptsize numerical values for various $b$}
& 
{\scriptsize 
 graphs of $w_{0,\Gamma},\;\tilde u_{0+}-w_{0,\Gamma},\;$ $b=10, 100, 1000, 10000$}
\end{longtable}\renewcommand{\arraystretch}{1}
 
The table shows that $\disp\lim_{b\to\infty}(\tilde u_{0+}-w_{0,\Gamma})=0$. Nevertheless, one may argue that, for low values of $b$, $w_{0,\Gamma}$ approximates loosely the exact solution $u_0$. Fortunately, this numerical analysis also suggests an easy way of obtaining truly sharp bounds {\em upper and lower} envelopes.

\begin{m-procedure}
Given a parameter $b$, one obtains sharp {\em upper and lower} envelopes by slightly increasing the value of the exponent $-10/3=-3.33\dots$. 
\end{m-procedure}
 For an upper bound, we replace $-3.33$ by a slightly greater $r$ ---let the resulting function be $w_{0,r}$---, such that the value at $x=1$ is positive; for lower bound, the value of $w_{0,r}$ at $x=0$ should be negative.
 
 \begin{m-example}
 (i) Let $b=10$.  A few trials yield: 
 \begin{itemize}
 \item for $r_-=-2.136804$, we have $w_{0,r_-}(1)=-\,6.2\cdot 10^{-9}<0$;
 \item for $r_+=-2.136805$, we have $w_{0,r_+}(1)=3.9\cdot 10^{-7}>0$. 
 \end{itemize}
 Thus, for $b=10$, the exact solution $u_0$ is squeezed between the envelopes $w_{0,r_-},w_{0,r_+}$, and the error is at most $4\cdot 10^{-7}$. 
 
\nit(ii) Let $b=10^6$. 
 \begin{itemize}
 \item $r_-=-3.129033$, $w_{0,r_-}=3.46\cdot 10^{-10}<0$;
 \item $r_+=-3.129034$, $w_{0,r_+}=8.16\cdot 10^{-10}>0$. 
 \end{itemize}
 For this choices, the exact solution is approximated with an error less than $1.2\cdot10^{-9}$. 
 
\nit(iii) Numerical data (see Table~\ref{tab:r}) suggests that the correct value of $r$ ---the one separating the upper and the lower envelopes--- follows the rule: 
\begin{m-eqn}{
r\approx -\frac{10}{3} + \frac{\rho}{\ln(b)},\quad\rho\approx 2.8.
}\label{eq:r}
\end{m-eqn} 
So the exponent $-10/3$ in equation~\eqref{eq:Gamma} defining $\Gamma$ is optimal, for arbitrary $b$.

\begin{longtable}{|l|l|l|l|}\caption{Precision of the R-K method, for $r$ given by~\eqref{eq:r}.}\label{tab:r}\\  \hline 
$b$ & $\begin{array}{c}\rho_-=2.84\\ w_{0,r_-}(1)\end{array}$   
& $\begin{array}{c}\rho_+=2.8\\ w_{0,r_+}(1)\end{array}$ & $\tilde u_{0+}(1)$ \\  \hline 
$10$          & -1.5 E-2      & -7.9 E-3       & 0.2087          \\  \hline 
$10^2$     & -8.6 E-5      & 1.4 E-3        & 0.0477          \\  \hline 
$10^3$     & -5.9 E-5      & 2.7 E-4        & 0.0103          \\  \hline 
$10^4$     & -2.2 E-5      & 4.9 E-5        & 2.2 E-3          \\  \hline 
$10^5$     & -6.0 E-6      & 9.5 E-6        & 4.8 E-4          \\  \hline 
$10^6$     & -1.4 E-6      & 1.8 E-6        & 1.0 E-4          \\  \hline 
$10^7$     & -3.4 E-7      & 3.7 E-7        & 2.2 E-5          \\  \hline 
$10^8$     & -7.9 E-8      & 7.6 E-8        & 4.8 E-6          \\  \hline 
$10^9$     & -1.7 E-8      & 1.5 E-8        & 1.0 E-6          \\  \hline 
$10^{10}$ & -3.9 E-9   & 3.2 E-9         & 2.2 E-7          \\  \hline 
\end{longtable}
The exact solution $u_0$ is squeezed between $w_{0,r_-}$ and $w_{0,r_+}$ (for $b\ges 16$). The value $b=10$ is an exception, one should take $\rho_+=2.75$ rather than $2.8$. The table suggests that the error of the computation is approximately $b^{-1}$. 
 \end{m-example}


\section{Numerical analysis}\label{sct:num}

So far, we developed a theoretical framework. Here we confront our analysis with numerical evidence, to probe our techniques.  

\subsection{Algorithm}\label{ssct:maple}

The code below\footnote{It's available also at:\\   \url{https://drive.google.com/file/d/1sPbb11_43iZ6V_Q0vdWHw5h-LJkKo4Oz/view?usp=sharing}.}, written in MAPLE, organizes the various formulae, and shows the interdependence of the quantities introduced so far. It runs on a usual computer, and allows the interested reader to experiment with the own favourite parameters. 

Note that if, for instance, one is interested in a partial lower bound for $y\in[T^{3/4},T]$, it suffices to change $sss$ ---first line of the code--- to $t^{3/4}$. The estimates will be more accurate. The value of $B$ is supposed to be at least $1$. 

\renewcommand{\arraystretch}{1.05}
 \begin{longtable}{l} 
 \caption{Code for upper/lower and global/partial envelopes}\\[-1ex] 
 \tiny $\text{restart}: N:=20000: Digits := N: \textbf{b := 70: t := 0.1: sss := $t^{1/2}$: } B := bt^{1.5}/\sqrt{5}: T := 1/t: L := 1.5\sqrt{2}B:$ 
\\ \tiny 
$BB := bb*tt^{1.5}/\sqrt{5}: TT := 1/tt: tL := L/qq: c0:=\text{write 30 iterations of the $\tanh(LL*)$ function}:$
\\ \tiny 
$ccm := \tanh(LL*c0+\tanh^{-1}(c0/zz)): ccp := tanh(LL*c0+\tanh^{-1}(1/zz)/(1-LL/\cosh(LL*c0)^2)): $
\\ \tiny 
$C := 2*cc/(1+ee-(1-ee)*cc): ZZ := (T^{3/2}+((1-ee)*(1/2))*C)/(1+((1-ee)*(1/2))*C):$
\\ \tiny 
$RR := ((1-(1/4)*q^2*(1+ee)^2)*5)*y^4+7*CC*(1-ee)*y^{5/2}-2*CC^2*(-ee^2+4*ee-1)*y$
\\ \tiny\qquad 
$-CC^3*ee*(1-ee)/y^(1/2)-CC^4*ee^2/y^2+q^2*(5*(1+ee)^2*(1/4)):$
\\ \tiny 
$nom:=5*y^{13/2}*q^2+16*CC^2*y^{7/2}-5*q^2*y^{5/2}+14*CC*y^5+2*CC^3*y^2 $   
\\ \tiny\qquad  
$-2*(70*CC*y^{23/2}*q^2+10*CC^3*y^{17/2}*q^2
+36*CC^3*y^{17/2}-70*CC*y^{15/2}*q^2 +36*CC^5*y^{11/2}$
\\ \tiny\qquad 
$+25*y^{13}*q^2 -10*CC^3*y^{9/2}*q^2 +60*CC^2*y^{10}*q^2-5*CC^4*y^7*q^2 +4*CC^3*y^{9/2}*r+9*CC^2*y^{10} $
\\ \tiny\qquad 
$-5*y^9*q^2*r+54*CC^4*y^7 -25*y^9*q^2+9*CC^6*y^4-60*CC^2*y^6*q^2+5*CC^4*y^3*q^2$
\\ \tiny\qquad $+8*CC^2*y^6*r-4*CC^4*y^3*r+5*y^5*q^2*r)^{1/2}:$
\\ \tiny
$ denom:=5*y^{13/2}*q^2-8*CC^2*y^{7/2}-5*q^2*y^{5/2}+4*CC^4*y^{1/2}-4*CC^3*y^2:$
\\ \tiny 
$ey :=-nom/denom:$ \# This is the root close to $1$ of the equation $R(\dots)=r$: 
\\ \tiny 
$mxder := (1-5*ss^4+(4+Dta)*ss^5)^{1/2}:qqss := 2*(1+(1-ee)*CC*ss^{1.5}-ee*CC^2*ss^3)/((1+ee)*mxder):$
\\ \tiny 
$maxders := subs(ss = sss, mxder):maxdert := subs(ss = t, mxder): minder := (1-5*t^4+4*t^5)^{1/2}:$
\\ \tiny
$qmin := 1-(3/5)*t^3: emax := 1-(4/5)*t^3+((8-3*t)*(1/5))*t^4:emin := 1-(4/5)*t^3+(3/5)*t^4:$
\\ \tiny\#\bf Define the global upper bound $typ$.
\\ \tiny
$tqp:=(((1-t^3)*(t^3+5))/(5*(1-t^4)))^{1/2}: tBp := B/tqp: typ := (\tanh(\tanh^{-1}(t^{1.5})+1.5*\sqrt{2}*tBp*(1-x)))^{-2/3}:$ 
\\ \tiny\#\bf Define the global lower bound $tym$.\\ \tiny 
$emm := 0.73: cmin := subs(LL = 0.99*L, c0): C0mm := 1/emm:  Dta := \frac{((subs(x =0 .25, typ)-1)}{(0.25*\sqrt{2}*B))^2}:$
\\ \tiny
$qqq := subs(ee = emm, CC = C0mm, ss = t, qqss): qmm := min(1.1, qqq): Lmm := subs(qq = qmm, tL): $
\\ \tiny 
$zmm := subs(ee = emm, cc = cmin, ZZ): cmm := evalf[N](subs(LL = Lmm, zz = zmm, ccm)):$
\\ \tiny
$Cmm := evalf[N](subs(ee = emm, cc = cmm, C)):$ 
\\ \tiny
\scalebox{.92}{$\disp tym := \Bigl(Cmm*\frac{T^{3/2}*(1+emm*\exp(2*Lmm*Cmm*(1-x)))+emm*Cmm*(\exp(2*Lmm*Cmm*(1-x))-1)} {T^{3/2}*(\exp(2*Lmm*Cmm*(1-x))-1)+Cmm*(\exp(2*Lmm*Cmm*(1-x))+emm)}\Bigr)^{2/3}:$}
\\ \tiny \#\bf  Define the partial upper bound $yp$.\\ \tiny 
$\disp qqp := \frac{2*(1+(1-ee)*CC*t^{1.5}-ee*CC^2*t^3)}{(1+ee)*minder}: e1p := subs(y = T, q = qmin, CC = cmin, r = 0, ey):$
\\ \tiny
$ q1p := subs(ee = emin, CC = cmin, qqp): L1p := subs(qq = q1p, tL): z1p := subs(ee = e1p, cc = 1, ZZ):$ 
\\ \tiny
$c1p := evalf[N](subs(LL = L1p/e1p, zz = z1p, ccp)): C1p := evalf[N](subs(ee = e1p, cc = c1p, C)):$
\\  \tiny
$R1 := evalf[N](subs(y = T, ee = e1p, q = q1p, CC = C1p, RR)): e2p := e1p: q2p := q1p+(1/15)*R1*t^4:$  
\\ \tiny
$L2p := subs(qq = q2p, tL): z2p := subs(ee = e2p, cc = 1, ZZ): c2p := evalf[N](subs(LL = L2p/e2p, zz = z2p, ccp)): $ 
\\ \tiny
$C2p := evalf[N](subs(ee = e2p, cc = c2p, C)):R2 := evalf[N](subs(y = T, ee = e2p, q = q2p, CC = C2p, RR)):  $
\\ \tiny
$e3p := e2p: q3p := q2p+(1/15)*R2*t^4:L3p := subs(qq = q3p, tL): z3p := subs(ee = e3p, cc = 1, ZZ): $
\\ \tiny
$c3p := evalf[N](subs(LL = L3p/e3p, zz = z3p, ccp)): C3p := evalf[N](subs(ee = e3p, cc = c3p, C)): $
\\ \tiny
$ep := e3p: qp := q3p: Cp := C3p: Lp := L3p: $
\\ \tiny 
$\disp yp := \Bigl(Cp*\frac{(T^{3/2}*(1+ep*\exp(2*Lp*Cp*(1-x)))+ep*Cp*(\exp(2*Lp*Cp*(1-x))-1))}
{(T^{3/2}*(\exp(2*Lp*Cp*(1-x))-1)+Cp*(\exp(2*Lp*Cp*(1-x))+ep))}\Bigr)^{2/3};$
\\ \tiny \#\bf  Define the partial lower bound $yp$.\\ \tiny 
$qmax := 1.25: ccmin := subs(LL = L/qmax, c0): CCmax := 1/emin: CCmin := ccmin: qtm := subs(ss = t, qqss):$ 
\\ \tiny 
$ qsm := subs(ss = sss, qqss): ee0m := emax: qq0m := subs(ee = emax, CC = CCmax, qsm):$ 
\\ \tiny 
$LL0m := subs(qq = qq0m, tL):zz0m := subs(ee = ee0m, cc = ccmin, ZZ):  $ 
\\ \tiny 
$cc0m := evalf[N](subs(LL = LL0m, zz = zz0m, ccm)):CC0m := evalf[N](subs(ee = ee0m, cc = cc0m, C)): $ 
\\ \tiny 
$ee1m := emax: qq1m := subs(ee = emax, CC = CC0m, qsm): LL1m := subs(qq = qq1m, tL): $
\\ \tiny 
$zz1m := subs(ee = ee1m, cc = ccmin, ZZ): cc1m := evalf[N](subs(LL = LL1m, zz = zz1m, ccm)): $ 
\\ \tiny 
$CC1m := evalf[N](subs(ee = ee1m, cc = cc1m, C)): 
em := ee1m: qqmm := qq1m: Cm := CC1m: Lm := LL1m:$
\\ \tiny 
$\disp ym := \Bigl(Cm*\frac{T^{3/2}*(1+em*\exp(2*Lm*Cm*(1-x)))+em*Cm*(\exp(2*Lm*Cm*(1-x))-1)}
{T^{3/2}*(\exp(2*Lm*Cm*(1-x))-1)+Cm*(\exp(2*Lm*Cm*(1-x))+em)}\Bigr)^{2/3}:$
\\ \tiny 
$plot([typ-tym,yp-ym],x=0..1);plot([yp,ym],x=0..1);$
\renewcommand{\arraystretch}{1} \end{longtable}\smallskip


\subsection{Numerical data}\label{ssct:num}

We are going to illustrate the enhanced precision of the partial bounds, compared to the global ones. Below we plotted the graphs of  various differences between the global and the partial envelopes, and also the plots of the envelopes themselves. The abbreviation B-L stands for `boundary layer'.\smallskip 

\begin{longtable}{rcr} 
\caption{Comparison of global and partial envelopes}\\ 
\hline 
\scriptsize $\tilde y_+-\tilde y_-$ (red),  $y_+-y_-$ (blue)
&
\scriptsize 
$\tilde y_+-y_+$ (red), $y_--\tilde y_-$ (blue)
&
\scriptsize 
$\tilde y_+$ (red), $\tilde y_-$ (blue)
\\ \hline && \\ 
\scriptsize 
$\begin{array}{l}b=70\\ t=0.1\\ B=0.99\\ bt^{0.25}=39\\ \text{B-L no} 
\end{array}$
\begin{minipage}[c]{0.245\textwidth}
\centering\includegraphics[width=.85\textwidth]{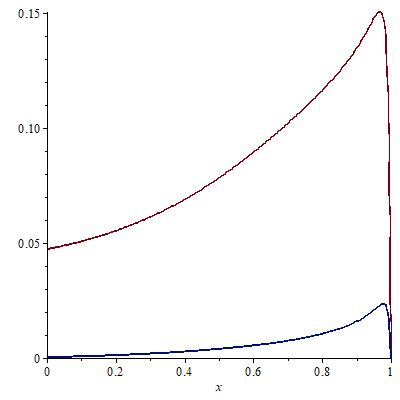}
\end{minipage}
&
\begin{minipage}[c]{0.245\textwidth}
\centering\includegraphics[width=.85\linewidth]{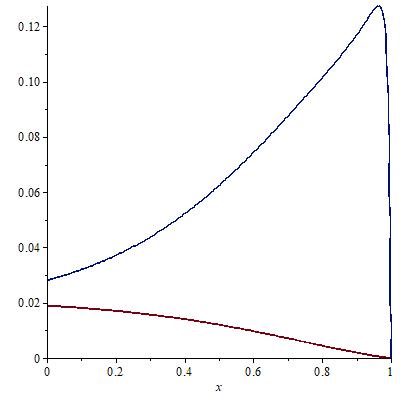}
\end{minipage}
&
\begin{minipage}[l]{0.245\textwidth}
\centering \includegraphics[width=.85\linewidth]{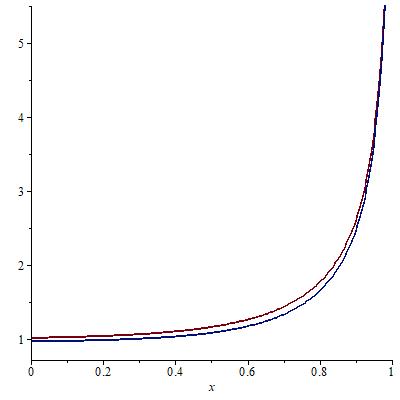}
\end{minipage}
\\ \hline && \\ 
\scriptsize$\begin{array}{l}b=500\\ t=0.1\\ B=7.07\\ bt^{0.25}=281\\ \text{B-L yes} \end{array}$
\begin{minipage}[c]{0.245\textwidth}
\centering 
\includegraphics[width=.85\linewidth]{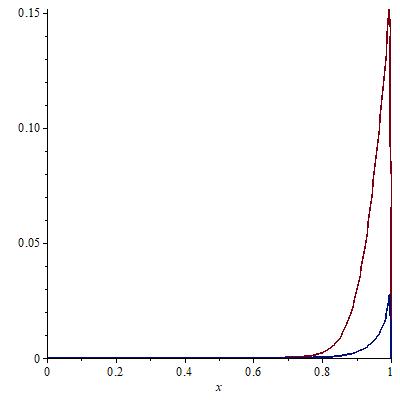}
\end{minipage}
&
\begin{minipage}[c]{0.245\textwidth}
\includegraphics[width=.85\linewidth]{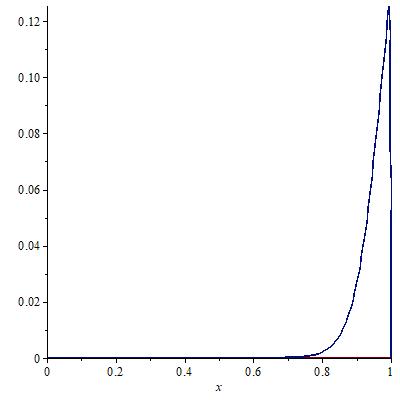}
\end{minipage}
&
\begin{minipage}[c]{0.245\textwidth}
\includegraphics[width=.85\linewidth]{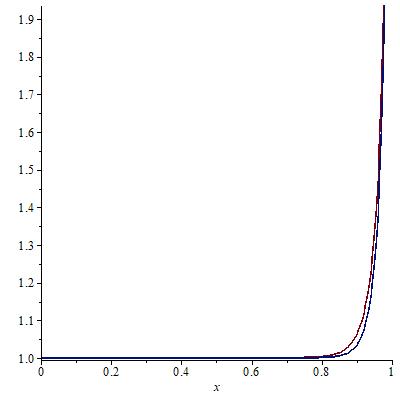}
\end{minipage} 
\\ \hline &&\\ 
\scriptsize $\begin{array}{l}b=10^6\\ t=2\cdot10^{-4}\\ B=1.26\\ bt^{0.25}=1.89\cdot10^5\\ \text{B-L yes} \end{array}$
\begin{minipage}[c]{0.245\textwidth}
\centering 
\includegraphics[width=.85\linewidth]{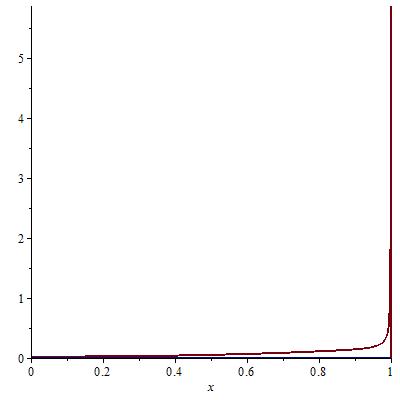}
\end{minipage}
&
\begin{minipage}[c]{0.245\textwidth}
\includegraphics[width=.85\linewidth]{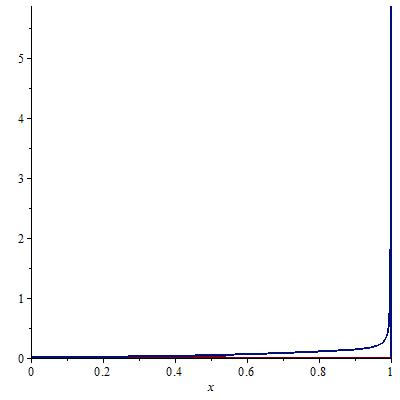}
\end{minipage}
&
\begin{minipage}[c]{0.245\textwidth}
\includegraphics[width=.85\linewidth]{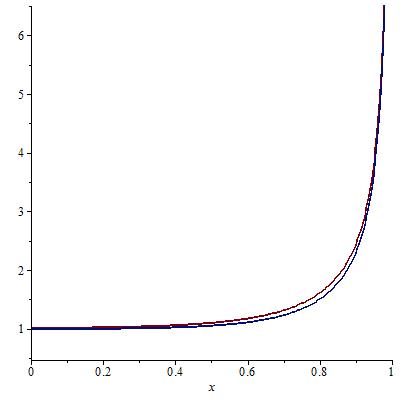}
\end{minipage} 
\\ \hline 
\end{longtable}\renewcommand{\arraystretch}{1}\medskip 

It's apparent that $\tilde y_+$ and $y_+$ are very close to each other ---in the second column of plots, the red graph basically are on the X-axis---, while $y_--\tilde y_-$ is not negligible; $y_-$ considerably improves $\tilde y_-$. Moreover, the difference $y_+-y_-$ of the partial envelopes ---the blue graphs in the first column--- is typically small. Thus the work done before pays off. We included the example, $b=10^6,t=2\cdot10^{-4}$, to show that our methods is functional even in numerically challenging situations.

The following table summarizes the computational error, the maximal difference between the upper and lower envelopes, for various values of the parameters $b,t$. The data suggests that the precision of the computation is related to the value $bt^{0.25}$ of the boundary layer: 
\begin{itemize}
\item When computing with global envelopes, the error (for $y$) increases with it, but the error for $u$ decreases.
\item When using partial envelopes, the precision increases (error decreases) with the value of the boundary layer, both for $y$ and $u$. The numerical data suggests that the number of exact decimals is of the same order as the boundary-value $bt^{0.25}$. This is precisely the desired behaviour.
\end{itemize}

\begin{table}[!ht]\renewcommand{\arraystretch}{1.2}
\caption{Maximal computational error for $y_\het$ and $u_\het$.}
\begin{tabular}{|rrrr|rr|rr|}
\hline
\multicolumn{4}{|l|}{parameters and errors}     & \multicolumn{2}{l|}{max. difference}      & \multicolumn{2}{l|}{max difference}          \\
\multicolumn{4}{|l|}{for $y$ and $u$}      & \multicolumn{2}{l|}{global envelopes}                 & \multicolumn{2}{l|}{partial envelopes}      \\ \hline
\multicolumn{1}{|c|}{$b$}       & \multicolumn{1}{c|}{$t$}      & \multicolumn{1}{c|}{$B$}    & \multicolumn{1}{c|}{$bt^{0.25}$}  
& \multicolumn{1}{c|}{\phantom{X}of $y$\phantom{X}} & \multicolumn{1}{c|}{of $u$} & \multicolumn{1}{c|}{of $y$}    & \multicolumn{1}{c|}{of $u$} \\ \hline
\multicolumn{1}{|r|}{500}     & \multicolumn{1}{r|}{0.1}    & \multicolumn{1}{r|}{7.07} & 281                     & \multicolumn{1}{r|}{0.16} & 0.016                     & \multicolumn{1}{r|}{0.027}   & 0.0027                    \\ \hline
\multicolumn{1}{|r|}{700}     & \multicolumn{1}{r|}{0.2}    & \multicolumn{1}{r|}{28}   & 468                     & \multicolumn{1}{r|}{0.11} & 0.022                     & \multicolumn{1}{r|}{0.024}   & 0.0048                    \\ \hline
\multicolumn{1}{|r|}{5000}    & \multicolumn{1}{r|}{0.01}   & \multicolumn{1}{r|}{2.23} & 1581                    & \multicolumn{1}{r|}{1.1}  & 0.011                     & \multicolumn{1}{r|}{0.017}   & 1.7 E-4                     \\ \hline
\multicolumn{1}{|r|}{10000}   & \multicolumn{1}{r|}{0.005}  & \multicolumn{1}{r|}{1.58} & 2659                    & \multicolumn{1}{r|}{2.17} & 0.01                      & \multicolumn{1}{r|}{0.013}  & 6.5 E-5                   \\ \hline
\multicolumn{1}{|r|}{1000000} & \multicolumn{1}{r|}{0.0002} & \multicolumn{1}{r|}{1.26} & 1.19 E+5          & \multicolumn{1}{r|}{5.87}  & 0.005                     & \multicolumn{1}{r|}{1.4 E-3} & 2.8 E-7                  
 \\ \hline  
\end{tabular}
\label{tab:maxerr-yu}
\end{table}\renewcommand{\arraystretch}{1}

Here we listed the worst possible error. But, for a specific value $x$, the computational error for $u_\het(x)$ is $t\cdot(y_+(1-x)-y_-(1-x))$, which can be much less than the maximal error, see Table~\ref{tab:precision}. Moreover, one can combine the envelope technique with the Runge-Kutta method in Section~\ref{ssct:shoot}, whose precision is $O(t\cdot\exp(-3B))$.

\begin{table}[!ht]\label{tab:precision}\renewcommand{\arraystretch}{1.1}
\caption{Upper and lower bounds for $u_\het$, for $b=500, t=0.1$.}
\begin{tabular}{lcc|c|}
\cline{4-4}  &  &   & number of   \\ \cline{1-3}
\multicolumn{1}{|c|}{$x$}   & \multicolumn{1}{c|}{$t\cdot y_+(x)$}     
& $t\cdot y_-(x)$       & exact decimals \\ \hline
\multicolumn{1}{|c|}{$10^{-2}$}  & \multicolumn{1}{c|}{0.314022890404343} & 0.311867729652350 & 2    \\ \hline
\multicolumn{1}{|c|}{$10^{-4}$}  & \multicolumn{1}{c|}{0.969598013211494} & 0.969224267767387 & 3    \\ \hline 
\multicolumn{1}{|c|}{$10^{-6}$}  & \multicolumn{1}{c|}{0.999684111921659} & 0.999680075693138 & 5  
 \\ \hline
\multicolumn{1}{|c|}{$10^{-8}$}  & \multicolumn{1}{c|}{0.999996839882428} & 0.999996799488355 & 6           
 \\ \hline
\multicolumn{1}{|c|}{$10^{-10}$} & \multicolumn{1}{c|}{0.999999968398700} & 0.999999967994756 & 8     \\  \hline  
\end{tabular}\renewcommand{\arraystretch}{1}
\end{table}

\section{Conclusion}

The article investigates the equation of heat conduction with radiation: 
$$
u''_\het=\bb^2(u^4_\het-\tu^4),\quad u_\het(0)=1,\;u_\het(1)=\tu\in(0,1).
$$
It combines analytical with numerical techniques. 

On the analytic side, we determined global and partial, upper and lower bounds for the exact solution $u_\het$. Also, we clarified the law, the value, governing the boundary layer property of $u_\het$. Moreover, we estimated the initial derivative $u'_\het(1)$ of the solution to the heat equation, which is essential for running the Runge-Kutta algorithm. The precision of the estimate for $u_\het$ is $O(t\cdot e^{-3B}).$

On the numerical side, we implemented the analytical arguments into the MAPLE computer program and produced numerical data for various parameters. The precision in estimating $u_\het$ ---the number of exact decimals--- has the same order of magnitude as the boundary-value $bt^{0.25}$.



\end{document}